\documentclass[11pt,a4paper]{article}
\usepackage[utf8]{inputenc}
\usepackage{amsmath}
\usepackage{amsfonts}
\usepackage{amssymb}
\usepackage{makeidx}
\usepackage{graphicx}
\usepackage{mathabx} 
\usepackage{dsfont} 
\usepackage[left=2cm,right=2cm,top=2cm,bottom=2cm]{geometry}
\usepackage{color}
\usepackage{verbatim} 
\usepackage{tocbibind} 
\usepackage{hyperref} 
\usepackage{mathtools} 
\usepackage{amsthm}
\usepackage{cancel}

\newtheorem{satz}{Lemma}[section]
\newtheorem{theorem}{Theorem}[section]
\newtheorem{definition}{Definition}[section]
\newtheorem{remark}{Remark}[section]
\newtheorem{proposition}{Proposition}[section]

\begin{document}

\setlength{\parindent}{0pt} 

\ \vspace{3cm}

\begin{center}
\Large\textbf{FLUID-RIGID BODY INTERACTION IN AN INCOMPRESSIBLE ELECTRICALLY CONDUCTING FLUID} \\[10mm]
\large{BARBORA BENE\v{S}OVÁ}$^1$, \large{\v{S}ÁRKA NE\v{C}ASOVÁ}$^2$, \large{JAN SCHERZ}$^{1,2,3}$,\\ \large{ANJA SCHLÖMERKEMPER}$^3$ \\[10mm]
\end{center}

\begin{itemize}
\item[$^1$] Department of mathematical analysis, Faculty of mathematics and physics, Charles University in Prague, Sokolovská 83, Prague 8, 18675, Czech Republic
\item[$^2$] Mathematical Institute, Academy of Sciences, \v{Z}itná 25, Prague 1, 11567, Czech Republic
\item[$^3$] Institute of Mathematics, University of Würzburg, Emil-Fischer-Str. 40, 97074 Würzburg, Germany
\end{itemize}      

\bigskip

\begin{center}
\Large\textbf{Abstract} \\[4mm]
\end{center}
      
We analyze a mathematical model that describes the interaction between an insulating rigid body and an incompressible electrically conducting fluid surrounding it. The model as well as the mathematical analysis involve the fields of fluid-structure interaction and of magnetohydrodynamics. Our main result shows the existence of a weak solution to the corresponding system of partial differential equations. The proof relies on a use of a time discretization via the Rothe method to decouple the system. This allows us to deal with test functions, depending on the position of the moving body and therefore on the solution of the system, in the weak formulation of the induction equation. The proof moreover makes use of the Brinkman penalization in order to cope with the mechanical part of the problem.

{\centering \section{Introduction}\par
}


We study a system of partial differential equations describing the movement of an insulating rigid body through an electrically conducting incompressible fluid. We prove the existence of a weak solution to this system. The insulating solid interacts mechanically with the fluid, hence the studied problem falls into the broad class of {\it fluid-structure interaction}. As the electrically conducting fluid further interacts with electromagnetic fields, it also constitutes a problem of \emph{magnetohydrodynamics}. While a number of mathematical works can be found in both those classes, a combination of the two seems to be missing. Possible applications include the interaction of both extra- and intracellular fluids with cell membranes in an organism. Even though membranes of cells are rather deformable, the study of the rigid body case can serve as a first step towards understanding the real-world situation.

\emph{Fluid-structure interaction} describes any interaction between a moving fluid and a rigid or deformable solid. In our case, we deal with the well-studied situation of a rigid body moving inside of a {\it viscous incompressible} fluid. For an introduction to the problem of a fluid coupled with a rigid body see \cite{G2,SER3}. First results on the existence of weak solutions until the first collision go back to the works of Conca, Starovoitov and Tucsnak \cite{CST}, Desjardins and Esteban \cite{desjardinsesteban}, Gunzburger, Lee and Seregin \cite{GLSE}, Hoffman and Starovoitov \cite{HOST}. The possibility of collisions has been addressed for example in \cite{tucsnak}, where the global-in-time existence of a weak solution in two dimensions is shown. Global-in-time existence results are also available in three dimensions for both incompressible and compressible fluids, see e.g.\ \cite{incompressiblefeireisl} and \cite{feireisl} respectively. Finally, we mention existence results in the setting of the Navier-slip boundary condition, c.f.\ \cite{compressiblecase}, and existence results on strong solutions, see e.g.\ \cite{GGH13,T,Wa}.

\textit{Magnetohydrodynamics} stands for the interaction of electrically conducting fluids with electromagnetic fields, see e.g.\ \cite{cabannes,kl}. This interplay is described mathematically by a coupling of the fluid equations with the Maxwell system \cite{landaulifshitz}. In this coupling, the resulting equations are further simplified subject to certain physical assumptions, which is referred to as the \textit{magnetohydrodynamic approximation}, c.f.\ \cite{davidson,eringenmaugin}. Results giving the existence of weak solutions to the magnetohydrodynamic model for incompressible and compressible fluids can be found in \cite{gbl} and \cite{sart} respectively. In \cite{blancducomet}, in addition to the electric conductivity, the fluid is assumed to be thermally conductive. The combination of an insulating rigid but also immovable object with an electrically conducting incompressible fluid has been considered in \cite{guermondminev2d} in two and in \cite{guermondminev} in three dimensions.

The model considered in the present work (see Section \ref{modeldescription}) is an extension of the model studied in the latter two articles. The novelty lies in the choice of the solid not as a fixed but as a freely moving rigid object, which causes various mathematical problems outlined below. We thus investigate the setting of a rigid non-conducting body moving in an incompressible electrically conducting fluid. We note that such a solid is trespassed by electromagnetic fields but not directly influenced by them. Thus, it may be viewed as vacuum from the electromagnetic point of view. Nonetheless, the electromagnetic fields influence indirectly its movement via the motion of the surrounding fluid. The setting we study can serve as a basis towards the study of more sophisticated systems involving for example a compressible fluid, an electrically conducting magnetic body or different types of boundary conditions.

Our main result yields the existence of weak solutions to the aforementioned model. In the weak formulation of the system we consider test functions which depend on the solid domain, i.e.\ on the solution of the system itself. While in pure fluid-structure interaction problems such test functions are standard and the resulting difficulties well-investigated, this is not the case for magnetohydrodynamical problems in moving domains. In our specific scenario, for the mechanical part we can rely on the well-studied Brinkman penalization method (see \cite{cottetmaitre}). However, for the problem resulting from the solution-dependent test functions in the induction equation, c.f.\ (\ref{633}) below, no such penalization method appears to be available. In order to deal with this problem we thus decided to decouple the system by discretizing it in the time variable via the Rothe method, c.f.\ \cite[Section 8]{roubicek}, and regularize it. This procedure constitutes the main novelty in our proof. Moreover, our proof adopts various methods from \cite{feireisl}. A detailed description of the idea of our proof is given in Section \ref{section3}.

\bigskip

The outline of the article is the following: In Section \ref{modeldescription} the model is described. The corresponding weak formulation together with the main result follows in Section \ref{section2}. The proof of the main result extends across Sections \ref{section3}--\ref{projectionlimit}: In Section \ref{section3} the approximation to the original system is presented, to which the existence of a solution is shown in Section \ref{existenceapproximation}. Sections \ref{rothelimit}--\ref{projectionlimit} deal with the limit passages required to return to the original equations. Finally, in Appendix~\ref{appendix} some auxiliary results are discussed. \\

\subsection{Model description} \label{modeldescription}

\begin{figure}[h!]
\begin{center}
\includegraphics[scale=0.5]{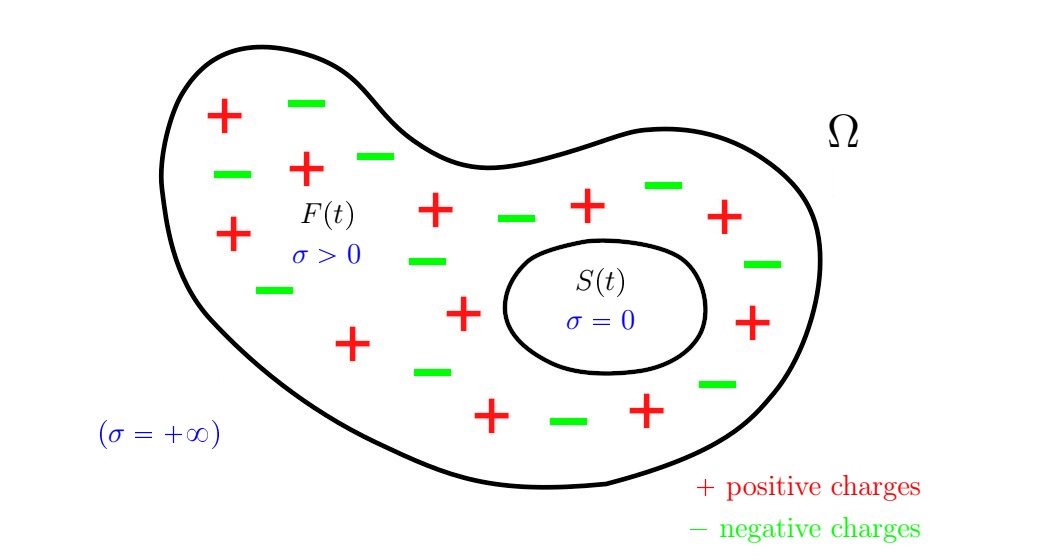} 
\caption{An insulating rigid body with domain $S(t)$ inside of an incompressible electrically conducting fluid with domain $F(t)$.} \label{figure}
\end{center}
\end{figure}

Let $\Omega \subset \mathbb{R}^3$ be a bounded domain occupied by a viscous nonhomogeneous incompressible fluid and a rigid body, let $T>0$ and set $Q:= (0,T) \times \Omega$. We denote the initial position of the body by $S=S(0)\subset \Omega$ and write $S(t)\subset \Omega$ for its position at any time $t \in [0,T]$, the movement of which can be expressed via an orientation preserving isometry. By $F(t) := \Omega \setminus \overline{S}(t)$ we denote the domain filled with the fluid at time $t$, see Figure \ref{figure}. Correspondingly, we divide $Q$ into the solid time-space domain
\begin{align}
Q^S := \left\lbrace (t,x) \in Q:\ x \in S(t) \right\rbrace. \nonumber
\end{align}

and its fluid counterpart $Q^F := Q \setminus \overline{Q}^S$, where $\overline{Q}^S$ denotes the closure of $Q^S$. We further split also any function $f$ defined on $Q$ into its fluid part $f^F$ and its solid part $f^S$,
\begin{align}
f(t,x) = \left\{
\begin{matrix}
f^F(t,x) \quad \quad \text{for } (t,x) \in Q^F \\
f^S(t,x) \quad \quad \text{for } (t,x) \in Q^S\ 
\end{matrix},
\right. \nonumber
\end{align}

whenever it is necessary to stress the difference. The motion of both the fluid and the body is then described by the velocity field $u:Q \rightarrow \mathbb{R}^3$, the density $\rho:Q \rightarrow \mathbb{R}$ and, in case of the fluid, also by the pressure $p = p^F:Q^F \rightarrow \mathbb{R}$. The electromagnetic effects in the system are characterized in by the magnetic induction $B:Q \rightarrow \mathbb{R}^3$, the magnetic field $H:Q \rightarrow \mathbb{R}^3$, the electric field $E:Q \rightarrow \mathbb{R}^3$ and the electric current $j:Q \rightarrow \mathbb{R}^3$. The evolution of the system is described by the following equations:
\begin{align}
\text{curl} H = j + J \quad \quad &\text{in } Q^F, \label{93} \\
\text{curl} H = 0 \quad \quad &\text{in } Q^S, \label{94} \\
\partial_t B + \text{curl} E = 0 \quad \quad &\text{in } Q, \label{95} \\
\text{div} E = 0 \quad \quad &\text{in } Q^S, \label{129} \\
\text{div} B = 0 \quad \quad &\text{in } Q,\label{876} \\
\nabla \cdot u = 0, \quad \partial _t \rho + u \cdot \nabla \rho = 0 \quad \quad &\text{in } Q^F \label{526} \\
\partial _{t} (\rho u) + \text{div} (\rho u \otimes u) + \nabla p = 2\nu \text{div} \mathbb{D} \left( u \right) + \rho g + \frac{1}{\mu} \text{curl}B\times B\quad \quad &\text{in } Q^F, \label{527} \\
m \frac{d}{dt} V(t) = \frac{d}{dt} \int_{\overline{S}(t)} \rho u\ dx = \int_{\partial\overline{S}(t)} \left[2\nu \mathbb{D}\left( u \right) - p\text{Id}\right] \cdot \text{n}\ d\sigma + \int _{\overline{S}(t)} \rho g \ dx,& \quad t \in [0,T] \label{528} \\
\frac{d}{dt}\left( \mathbb{J}(t)w(t) \right) = \frac{d}{dt} \int _{\overline{S}(t)} \rho \left(x - X\right) \times u \ dx\ & \nonumber \\
= \int_{\partial \overline{S}(t)} (x - X) \times \left[2\nu \mathbb{D} \left( u \right) - p\text{Id} \right] \text{n}\ d\sigma + \int _{\overline{S}(t)} \rho \left(x - X\right) \times g\ dx,& \quad t \in [0,T] \label{529}
\end{align}
supplemented by the relations
\begin{align}
&j = \sigma (E + u \times B) \quad \text{in } Q,\quad \sigma = \left\{
                \begin{array}{ll}
                  \sigma ^F > 0 \ \ \ &\text{in } Q^F,\\
                  \sigma ^S = 0 &\text{in } Q^S,
                \end{array}
              \right. \label{877} \\
&\quad \quad \quad \quad \ \ B = \mu H, \quad \quad \mu > 0 \quad \text{in } Q \label{878}
\end{align}
and completed by the boundary and interface conditions
\begin{align}
B(t) \cdot \text{n} &= 0\quad \text{on } \partial \Omega,\quad \quad &&\ \ \ \ \ \ \ \ \ B^F(t) - B^S(t) = 0 \quad \text{on } \partial S(t), \label{831} \\
E(t) \times \text{n} &= 0 \quad \text{on } \partial \Omega,\quad \quad &&\left(E^F(t) - E^S(t)\right) \times \text{n} = 0 \quad \text{on } \partial S(t), \label{832} \\
u(t) &= 0\quad \text{on } \partial \Omega,\quad \quad &&\ \ \ \ \ \ \ \ \ \ u^F(t) - u^S(t) = 0 \quad \text{on } \partial S(t). \label{833}
\end{align}

The electromagnetic part (\ref{93})--(\ref{876}) of the system is the Maxwell system, simplified according to the magnetohydrodynamic approximation in $Q^F$ and adjusted to the non-conductivity of the solid in $Q^S$, c.f.\ \cite{guermondminev2d}, \cite{guermondminev}. The mechanical part (\ref{526})--(\ref{529}) consists of the balance of mass and momentum for the fluid and solid respectively, c.f.\ for example \cite{incompressiblefeireisl}. For the Maxwell system we remark that Ampère's law (\ref{93}) in the fluid domain contains, as in \cite{guermondminev2d} and \cite{guermondminev}, an additional source term $J$. Moreover, the Maxwell system is supplemented by Ohm's law (\ref{877}), which determines the effect of the fluid motion on the electromagnetic quantities via the electrical conductivity $\sigma$. For the constitutive relation (\ref{878}) (c.f.\ \cite[Section 5.8]{jackson}) we point out that in our case the magnetic permeability $\mu$ is chosen as a constant value in the whole domain $Q$. As the magnetic permeability depends on the material, this is in general not physically accurate but a simplification required for the continuity of $B$ across the interface stated in (\ref{831}). The latter condition is needed to ensure that $B$ is an element of some Sobolev space on $Q$, c.f.\ (\ref{776}) below. The remaining boundary and interface conditions stated in (\ref{831}) and (\ref{832}) are standard.

In the mechanical part of the above system, the continuity equation and incompressibility condition (\ref{526}) and the momentum equation (\ref{527}) constitute the incompressible Navier-Stokes system. Besides the given external force $g$, (\ref{527}) contains the reduced Lorentz force $\frac{1}{\mu} \text{curl}B\times B$ resulting from the electromagnetic interaction. Moreover, the operator $\mathbb{D}$ denotes the symmetric part of the gradient, $\mathbb{D}(u) := \frac{1}{2} \nabla u + \frac{1}{2} (\nabla u)^T$, and $\nu > 0$ is the viscosity coefficient. The balance of linear momentum (\ref{528}) and angular momentum (\ref{529}) determine the translational velocity $V$ and the rotational velocity $w$ of the rigid body respectively and hence its total velocity as
\begin{equation}
u(t,x) = V(t) + w(t) \times \left(x - X(t)\right) \quad \text{on } Q^S. \nonumber
\end{equation}
The further notation in (\ref{528}) and (\ref{529}),
\begin{align}
m :=& \int _{\overline{S}(t)} \rho (t,x)\ dx,\quad X(t):= \frac{1}{m}\int _{\overline{S}(t)} \rho (t,x)x\ dx, \nonumber \\
\mathbb{J}(t)a \cdot b :=& \int _{\overline{S}(t)}\rho (t,x)\left[ a \times \left( x - X(t) \right) \right] \cdot \left[ b \times \left( x - X(t) \right) \right]\ dx,\quad a,b \in \mathbb{R}^3, \nonumber
\end{align}
represent the mass, the center of mass and the inertia tensor of the rigid body respectively, while $\text{n}$ denotes the outer unit normal vector on $\partial \Omega$ and $\partial S(t)$. The right-hand sides of (\ref{528}) and (\ref{529}) show that the movement of the rigid body is driven by the volume force $g$ and the Cauchy stress $2\nu \mathbb{D}(u) - p\text{Id}.$ Moreover, the fluid-structure interaction is incorporated in the no-slip interface condition on $\partial S(t)$ in (\ref{833}). This, as well as the no-slip boundary condition on $\partial \Omega$ in (\ref{833}), constitutes a common assumption for the interaction between fluids and rigid bodies, c.f.\ for example \cite{incompressiblefeireisl}, \cite{tucsnak}.

\subsection{Weak formulation and main result} \label{section2}

In order to define a weak formulation and state the main result, we first introduce some more notation. We use the standard Lebesgue-, Sobolev- and Bochner spaces and in addition the spaces
\begin{align}
V^r(\Omega) := H_\text{div}^{r}\left(\Omega \right)\quad \text{for } r \geq 0,\quad \quad V^r_0(\Omega) := \left\lbrace v \in H_\text{div}^{r}\left(\Omega \right):\ v|_{\partial \Omega} = 0 \right\rbrace \quad \text{for } r > \frac{1}{2}, \nonumber
\end{align}

where $H_\text{div}^{r}(\Omega)$ denotes the space of functions in the (potentially fractional) Sobolev space $H^{r}(\Omega)$ which are in addition divergence-free. By $S$ we denote a subset of $\Omega$ such that
\begin{align}
S \ \text{is open, bounded and connected, } S \neq \emptyset,\ |\partial S|=0 \text{ and } \text{dist}\left( S, \partial \Omega \right)>0. \label{777}
\end{align}

In the following we will describe the position of the rigid body through the use of a characteristic function. More precisely, for a function $\chi: \mathbb{R}^3 \rightarrow \left\lbrace 0, 1 \right\rbrace$ we write
\begin{align}
S(\chi) := \left\lbrace x \in \mathbb{R}^3:\ \chi(x) = 1 \right\rbrace \nonumber
\end{align}

and if $\chi(t)=\chi(t,\cdot)$ is a characteristic function for all $t \in [0,T]$, we write
\begin{align}
Q^S(\chi, T) := \left\lbrace (t,x) \in [0,T] \times \mathbb{R}^3:\ \chi(t,x) = 1 \right\rbrace = \left\lbrace (t,x) \in [0,T] \times \mathbb{R}^3:\ x \in S(\chi(t)) \right\rbrace. \nonumber
\end{align}

We further introduce the space of test functions
\begin{align}
\mathcal{T}(\chi, T) :=& \left\lbrace \phi \in \mathcal{D}\left([0,T)\times \Omega \right):\ \text{div}\phi = 0,\ \mathbb{D}(\phi) = 0\ \ \ \text{on an open neighbourhood of } \overline{Q}^S(\chi, T)\right\rbrace, \nonumber
\end{align}

meaning that for any $\phi \in \mathcal{T}(\overline{Q}^S(\chi, T))$ there is $\sigma > 0$ such that
\begin{align}
\mathbb{D}(\phi) = 0 \quad \text{in } \left\lbrace (t,x) \in (0,T)\times \Omega :\ \text{dist}\left( (t,x), \overline{Q}^S(\chi, T) \right) < \sigma \right\rbrace. \label{933}
\end{align}

Similarly, we define
\begin{equation}
Y(\chi, T):= \left\lbrace b \in \mathcal{D}\left([0,T)\times \Omega \right): \text{curl}\ b= 0 \quad \text{on an open neighbourhood of } \overline{Q}^S(\chi, T)\right\rbrace. \label{936}
\end{equation}

We can now introduce the notion of weak solutions to (\ref{93})--(\ref{832}) in the following way:
\begin{definition}
\label{weaksolutions}
Let $T > 0$, let $B_0,u_0 \in L^2(\Omega;\mathbb{R}^3)$ with $\operatorname{div} u_0= \operatorname{div} B_0=0$ and let $\rho _0 \in L^\infty(\Omega;\mathbb{R})$. Let further $S \subset \Omega$ satisfy (\ref{777}) and let $\chi_0 \in L^\infty(\Omega; \mathbb{R})$ denote the characteristic function of $S$. Then the set of functions $\left\lbrace \chi, \rho, u, B\right\rbrace$ is called a weak solution of the problem (\ref{93})--(\ref{832}) on $[0,T]$, if
\begin{align}
&\chi \in C\left([0,T];L^p\left(\Omega; \left\lbrace 0,1 \right\rbrace \right)\right) \quad \forall 1 \leq p < \infty, \label{665} \\
&\rho \in C\left([0,T];L^p\left(\Omega; \mathbb{R}\right)\right) \quad \forall 1 \leq p < \infty, \label{635} \\
&u \in L^\infty \left(0,T;L^2\left(\Omega; \mathbb{R}^3\right)\right) \bigcap L^2(0,T;V^1_0(\Omega)), \quad \mathbb{D}(u) = 0 \ \text{a.e. in } Q^S(\chi, T), \label{631} \\
&B \in L^\infty \left(0,T;L^2\left(\Omega; \mathbb{R}^3\right)\right) \bigcap L^2(0,T;V^1(\Omega)),\quad \operatorname{curl} B = 0 \ \text{a.e. in } Q^S(\chi, T),\quad B \cdot \text{n}= 0 \ \text{on } \partial \Omega, \label{776}
\end{align}

and
\begin{align}
-\int_0^{T} \int_{\Omega} \chi \partial_t \Theta dxdt - \int _{\Omega } \chi _0 \Theta (0,x)\ dx =& \int_0^T \int_\Omega \left( \chi u \right) \cdot \nabla \Theta\ dxdt, \label{666} \\
-\int_0^{T} \int _{\Omega} \rho \partial_t \psi dxdt - \int_{\Omega} \rho_0 \psi(0,x )\ dx =& \int_0^{T}\int_\Omega (\rho u) \cdot \nabla \psi \ dxdt  , \label{636} \\
- \int_0^{T} \int_\Omega \rho u \cdot \partial_t \phi\ dxdt - \int_\Omega \rho _0 u_0 \cdot \phi(0,x)\ dx =& \int _0^{T} \int _\Omega \left(\rho u \otimes u \right): \nabla \phi - 2\nu \mathbb{D}(u):\nabla \phi \nonumber \\
&+ \rho g\cdot \phi + \frac{1}{\mu}\left( \operatorname{curl} B \times B \right) \cdot \phi \ dxdt, \label{632} \\
 - \int _0^{T} \int_\Omega B \cdot \partial_t b\ dxdt - \int_\Omega B_0 \cdot b(0,x)\ dx =& \int _0^{T} \int _\Omega \left[ -\frac{1}{\sigma \mu} \operatorname{curl} B + u \times B + \frac{1}{\sigma} J \right] \cdot \operatorname{curl} b\ dxdt \label{633}
\end{align}

for all $\Theta$, $\psi \in \mathcal{D}([0,T)\times \Omega)$, $\phi \in \mathcal{T}(\chi, T)$ and $b \in Y(\chi, T)$. 

\end{definition}

We can now state our main result:
\begin{theorem}
\label{mainresult}

Let $T>0$, $\Omega \subset \mathbb{R}^3$ be a simply connected bounded domain. Let further $\chi_0$, $\rho_0$, $u_0$, $B_0$ and $S$ be as in Definition \ref{weaksolutions} and assume $\partial \Omega$, $\partial S$ to be of class $C^2 \bigcap C^{0,1}$. Finally, assume $g,J \in L^\infty(Q;\mathbb{R}^3)$ and assume that $\sigma$, $\mu$, $\nu$, $\underline{\rho}$, $\overline{\rho}$ are positive constants with
\begin{align}
0 < \underline{\rho} \leq \rho_0 \leq \overline{\rho} < \infty \quad \text{a.e. on } \Omega. \nonumber
\end{align}

Then there exist $T'>0$ and a weak solution to the problem (\ref{93})--(\ref{832}) on $[0,T']$ in the sense of Definition \ref{weaksolutions}, satisfying the energy inequality
\begin{align}
&\int_\Omega \frac{1}{2} \rho (\tau) |u(\tau)|^2 + \frac{1}{2}|B(\tau)|^2\ dx + \int_0^\tau \int_\Omega 2\nu \left| \nabla u(t,x) \right|^2 + \frac{1}{\sigma \mu^2} \left| \operatorname{curl} B(t,x) \right|^2\ dxdt \nonumber \\
\leq& \int_\Omega \frac{1}{2} \rho_0 |u_0|^2 + \frac{1}{2}|B_0|^2\ dx + \int_0^\tau \int_\Omega \rho (t,x) g(t,x)\cdot u(t,x) + \frac{1}{\sigma}J(t,x) \cdot \operatorname{curl} B(t,x)\ dxdt \label{638}
\end{align}

for almost all $\tau \in [0,T']$. Moreover, there is an orientation preserving isometry $X(s;t,\cdot):\mathbb{R}^3 \to \mathbb{R}^3$ such that
\begin{align}
S\left(\chi (t)\right) = X\left(s;t, S\left(\chi (s)\right)\right) \quad \forall s,t \in [0,T'] \label{780}
\end{align}

and $T'$ can be chosen such that
\begin{align}
T' = \sup \Big\lbrace \tau \in [0,T]:\ \operatorname{dist}\left( S(\chi(t)),\partial \Omega \right)>0\quad \forall t \leq \tau \Big\rbrace. \label{779}
\end{align}

\end{theorem}

\begin{remark}
Due to the transport theorem by DiPerna and Lions \cite{dipernalions}, the solution $\rho \in L^2((0,T)\times \Omega)$ to the continuity equation (\ref{636}), given by Theorem \ref{mainresult}, also solves the renormalized continuity equation
\begin{align}
\partial_t \beta(\rho) + u \cdot \nabla \beta(\rho) = 0 \quad \text{in } \mathcal{D}'((0,T) \times \Omega) \label{925}
\end{align}
for any bounded $\beta \in C^1(\mathbb{R})$ vanishing near $0$ and such that also $(\beta '(1 + |\cdot|))^{-1}$ is bounded.
\end{remark}

The proof of Theorem \ref{mainresult} will be accomplished via an approximation method in Sections \ref{existenceapproximation}--\ref{projectionlimit} and is outlined in the following section.

{\centering \section{Approximate system}\par \label{section3}
}

We introduce the approximation to the original system. Then again the limit of the solution to this approximation is supposed to give us the solution to the original problem. Our approximation consists of three different levels, characterized through three parameters $\Delta t, \epsilon, \eta>0$:
\begin{itemize}
\item On the $\Delta t$-level, we have a time discretization by the Rothe method, c.f.\ \cite[Section 8.2]{roubicek}: For $\Delta t>0$ such that $\frac{T}{\Delta t}\in \mathbb{N}$ we split up the interval $[0,T]$ into the discrete times $k\Delta t$, $k=1,...,\frac{T}{\Delta t}$.
\item On the $\epsilon$-level, we add several regularization terms to the system, which help us to solve the approximate system and pass to the limit as $\Delta t \rightarrow 0$.
\item On the $\eta$-level we add a penalization term to the momentum equation, which guarantees us that after passing to the limit in $\eta \rightarrow 0$, the limit velocity will coincide - on the solid part of the domain - with the rigid velocity of the body.
\end{itemize}

We now present the full approximate system, containing all three approximation levels, and give a more detailed explanation afterwards: Assuming for some discrete time $k \Delta t$, $k \in \{1,...,\frac{T}{\Delta t}\}$, the solution at time $(k-1)\Delta t$, indexed by $k-1$, to be given and defining the test function space
\begin{align}
W^k(\Omega) :=& \bigg \lbrace b\in H^{2}\left(\Omega \right):\ b \cdot \text{n}|_{\partial \Omega} = 0,\ \text{curl}\ b = 0\ \text{in } \Big \lbrace x \in \Omega:\ \chi_{\Delta t}^k(x) = 1\Big \rbrace \bigg \rbrace, \label{937}
\end{align}

we seek functions
\begin{align}
\rho_{\Delta t}^k:\ \Omega \rightarrow \mathbb{R},\quad u_{\Delta t}^k, \ B_{\Delta t}^k:\ \Omega \rightarrow \mathbb{R}^3,\quad \chi_{\Delta t,k}:\ [(k-1)\Delta t,k\Delta t] \times \Omega \rightarrow \mathbb{R}, \nonumber
\end{align}

satisfying the discrete system at time $k \Delta t$,
\begin{align}
- \int_{(k-1)\Delta t}^{k\Delta t} \int_{\mathbb{R}^3 } \chi_{\Delta t,k} \partial_t \Theta \ dxdt =&\int _{\mathbb{R}^3 } \chi _{\Delta t}^{k-1} \Theta ((k-1) \Delta t,x)\ dx - \int_{\mathbb{R}^3} \chi^k_{\Delta t} \Theta (k \Delta t,x)\ dx \nonumber \\
&+\int_{(k-1)\Delta t}^{k \Delta t} \left(\chi_{\Delta t,k} \Pi_{\Delta t}^{k-1} \right) \cdot \nabla \Theta \ dxdt, \label{236} \\
-\int _\Omega \frac{\rho _{\Delta t}^k - \rho _{\Delta t}^{k-1}}{\Delta t}\psi  dx =& \int_\Omega u_{\Delta t}^{k-1} \cdot \nabla \rho _{\Delta t}^k \psi + \epsilon \nabla \rho^k_{\Delta t} \cdot \nabla \psi \ dx, \label{8} \\
-\int_\Omega \frac{\rho _{\Delta t}^k u_{\Delta t}^k - \rho _{\Delta t}^{k-1} u_{\Delta t}^{k-1}}{\Delta t} \cdot \phi \ dx =& \int_\Omega \left[ \operatorname{div}\left( \rho _{\Delta t}^{k} u_{\Delta t}^{k-1} \otimes u_{\Delta t}^k \right) - 2\nu \operatorname{div} \left( \mathbb{D}u_{\Delta t}^k \right) + \epsilon \nabla u_{\Delta t}^k \nabla \rho_{\Delta t}^k \right] \cdot \phi \ dx \nonumber \\
&+ \int_\Omega \epsilon \Delta u_{\Delta t}^k \Delta \phi + \left[ \frac{1}{\eta} \rho _{\Delta t}^{k-1} \chi_{\Delta t}^{k} \left( u_{\Delta t}^{k-1} -  \Pi_{\Delta t}^{k-1} \right) - \rho_{\Delta t}^{k-1} g_{\Delta t}^k \right] \cdot \phi\ dx \nonumber \\
&- \int_\Omega \frac{1}{\mu} \left( \operatorname{curl} B_{\Delta t}^{k-1} \times B_{\Delta t}^{k-1} \right) \cdot \phi \ dx, \label{57} \\
-\int_\Omega \frac{B_{\Delta t}^k - B_{\Delta t}^{k-1}}{\Delta t} \cdot b\ dx =& \int_\Omega \left[ \frac{1}{\sigma \mu} \operatorname{curl} B_{\Delta t}^k  - u_{\Delta t}^k \times B_{\Delta t}^{k-1} + \frac{\epsilon}{\mu ^2} \left| \operatorname{curl} B_{\Delta t}^k \right|^2\operatorname{curl} B_{\Delta t}^k \right] \cdot \operatorname{curl} b\ dx \nonumber \\
&+ \int_\Omega -\frac{1}{\sigma} J_{\Delta t}^k \cdot \operatorname{curl} b + \epsilon \operatorname{curl}\left(\operatorname{curl} B_{\Delta t}^k\right) \cdot \operatorname{curl}\left(\operatorname{curl} b\right)\ dx \label{116}
\end{align}

for all $\Theta \in \mathcal{D}([(k-1)\Delta t,k\Delta t]\times \mathbb{R}^3)$, $\psi \in H^{1}(\Omega)$, $\phi \in V_0^2(\Omega)$ and $b \in W^k(\Omega)$. Here, the functions $\chi^k_{\Delta t}$ and $\Pi_{\Delta t}^{k-1}$, introduced in the equations (\ref{236}) and (\ref{57}), are defined by:
\begin{align}
\chi_{\Delta t}^{k} := \chi_{\Delta t,k}(k\Delta t),\quad &\Pi _{\Delta t}^{k-1} = (u_G)^{k-1}_{\Delta t} + \omega ^{k-1}_{\Delta t} \times \left(x - a^{k-1}_{\Delta t} \right) \label{222}
\end{align}

and
\begin{align}
(u_G)_{\Delta t}^{k-1} &:= \frac{\int _{\mathbb{R}^3} \rho _{\Delta t}^k \chi_{\Delta t}^{k-1} u_{\Delta t}^{k-1}\ dx }{\int _{\mathbb{R}^3} \rho _{\Delta t}^k \chi_{\Delta t}^{k-1}\ dx}, \label{835} \\
\omega _{\Delta t}^{k-1} &:= \left( I_{\Delta t}^{k-1} \right)^{-1} \int _{\mathbb{R}^3} \rho _{\Delta t}^k \chi_{\Delta t}^{k-1} \left( x - a_{\Delta t}^{k-1} \right) \times u_{\Delta t}^{k-1}\ dx ,\label{857} \\
I_{\Delta t}^{k-1} &:= \int _{\mathbb{R}^3} \rho _{\Delta t}^k \chi_{\Delta t}^{k-1} \left( | x - a_{\Delta t}^{k-1} |^2\text{id}  - \left( x - a_{\Delta t}^{k-1} \right) \otimes \left( x - a_{\Delta t}^{k-1} \right) \right)\  dx,\label{836} \\
a_{\Delta t}^{k-1} &:= \frac{\int _{\mathbb{R}^3} \rho _{\Delta t}^k \chi_{\Delta t}^{k-1} x\ dx }{\int _{\mathbb{R}^3} \rho _{\Delta t}^k \chi_{\Delta t}^{k-1}\ dx}. \label{837}
\end{align}

In order to keep the latter terms well-defined, we extend the functions $\rho_{\Delta t}^l$ by $\underline{\rho}$ and $u_{\Delta t}^l$ by $0$ outside of $\Omega$ for any $l=0,...,k$. Moreover, since the given $L^\infty$-functions $g$ and $J$ are not necessarily defined in the discrete times, we regularize them as in \cite[(7.10)]{roubicek},
\begin{align}
g_\gamma (t):= \int_0^T \theta _\gamma \left( t + \xi_\gamma (t) - s \right)g(s)\ ds,\quad J_\gamma (t):= \int_0^T \theta _\gamma \left( t + \xi_\gamma (t) - s \right)J(s)\ ds,\quad \xi_\gamma (t) := \gamma \frac{T-2t}{T}. \nonumber
\end{align}

Here $\theta_\gamma :\mathbb{R}\to \mathbb{R}$ is a mollifier with support in $[-\gamma, \gamma]$ and $\xi_\gamma$ has the purpose of shifting the support of $s \mapsto \theta _\gamma \left( t - s \right)$ into $[0,T]$ for any fixed $t \in [0,T]$. Then we choose $\gamma = \gamma ( \Delta t),$ $\gamma (\Delta t)\rightarrow 0$ for $\Delta t \rightarrow 0$ and define the discrete approximations $g_{\Delta t}^k$ and $J_{\Delta t}^k$ of $g$ and $J$ from (\ref{57}) and (\ref{116}) by
\begin{align}
g_{\Delta t}^k := g_{\gamma (\Delta t)}(k\Delta t),\quad J_{\Delta t}^k := J_{\gamma (\Delta t)}(k\Delta t). \label{840}
\end{align}

The idea behind the time discretization is to decouple the system in order to solve the induction equation, in which the test functions depend on the position of the body and thus on the overall solution, c.f.\ (\ref{936}). Indeed, using time-lagging functions, at each fixed discrete time we can first determine the solid domain and subsequently choose the test functions at that specific time accordingly. The discrete induction equation can then be solved by standard methods.

The characteristic function $\chi_{\Delta t,k}$, which is immediately constructed as a time-dependent function, represents an exception in this discrete system. The reason for this is that for the determination of the solid domain it should take only the values $0$ and $1$. Inspired by \cite{gigli}, we can guarantee this by constructing $\chi_{\Delta t,k}$ by solving a classical transport equation on the small interval $[(k-1)\Delta t, k\Delta t]$, in case of a discrete transport equation we might lose the property.

Next we note that the mapping $\Pi_{\Delta t}^{k-1}$ is, by definition, a rigid velocity field. After the limit passage in the time discretization, we will see that the limit of $\overline{\Pi}_{\Delta t}'$ is actually the projection of the velocity onto a rigid velocity field. This comes into play in the penalization term from the $\eta$-level, namely the term
\begin{align}
\frac{1}{\eta} \rho _{\Delta t}^{k-1} \chi_{\Delta t}^{k} \left( u_{\Delta t}^{k-1} - \Pi_{\Delta t}^{k-1} \right), \nonumber
\end{align}

from (\ref{57}). As mentioned before, we can use this term to infer that after letting $\eta \rightarrow 0$ the limit velocity coincides, in the solid area, with the velocity of the rigid body, which is what we require to obtain (\ref{666}). This penalization method, which is known as Brinkman penalization, is discussed rigorously in \cite{cottetmaitre}. It can be considered as an extension of the penalty method used in \cite{abf} for a fluid-structure interaction problem in which the movement of the solid is prescribed. It further finds use in \cite{alternativerigidbodies}, where the solid is additionally deformable and self-propelled and it is moreover of interest for finite element approaches to the problem, c.f.\ \cite{cc}, \cite{hl}.

Finally, it remains to discuss the regularization terms from the $\epsilon$-level. In the continuity equation, the Laplacian of the density is added, which allows us to show an upper bound for $\rho$ as well as some bound away from $0$. This is needed because such a bound cannot be guaranteed from the discrete version of the transport equation. The quantity $\epsilon \nabla u_{\Delta t}^k \nabla \rho_{\Delta t}^k$ in the momentum equation compensates for this term in the energy inequality. The second new term in this equation, $\epsilon \Delta^2 u_{\Delta t}^k$, is needed for passing to the limit in $\Delta t \rightarrow 0$. Moreover, we have two regularization terms in the induction equation, the $4$-th curl of $B_{\Delta t}^k$ and the $4$-double curl, $\text{curl} (| \text{curl} B_{\Delta t}^k |^2 \text{curl}B_{\Delta t}^k)$. The first one is used for the construction of $B_{\Delta t}^k$ via a weakly continuous coercive operator, while we require the latter one in the energy inequality: since the mixed terms from the momentum and the induction equation are chosen from distinct discrete times, they do not cancel each other as in the continuous system. However, the $4$-double-curl enables us to absorb the problematic terms into the positive left-hand side, so that we can get the uniform bounds needed for the limit passage as $\Delta t \rightarrow 0$. We complement the equations by the relations
\begin{align}
\operatorname{div} u_{\Delta t}^k = \operatorname{div} B_{\Delta t}^k = 0\quad \text{in } \Omega,&\quad \quad \operatorname{curl} B_{\Delta t}^k = 0 \quad \text{in } \left\lbrace x\in \Omega:\ \chi_{\Delta t}^k(x) = 1\right\rbrace, \nonumber \\u_{\Delta t}^k = 0 \quad \text{on } &\partial \Omega,\quad \quad B_{\Delta t}^k \cdot \text{n} = 0 \quad \text{on } \partial \Omega, \label{893} \\
\rho_{\Delta t}^0 = \rho _0,\quad \quad \chi_{\Delta t,k}(0) &= \chi_0,\quad \quad u_{\Delta t}^0 = u_0,\quad \quad B_{\Delta t}^0 = B_0. \nonumber
\end{align}

{\centering \section{Existence of the approximate solution} \label{existenceapproximation} \par 
}

In this section we prove the existence of a solution to the approximate system. To this end, we first introduce another function space for fixed discrete time indices $k$:
\begin{align}
Y^k(\Omega) :=&\bigg \lbrace b \in L^2\left(\Omega; \mathbb{R}^3\right):\ b \cdot \text{n}|_{\partial \Omega} = 0,\ \text{div}\ b=0\ \text{in } \Omega,\ \text{curl}\ b = 0\ \text{in } \Big \lbrace x\in \Omega:\ \chi_{\Delta t}^k(x) = 1\Big \rbrace , \nonumber \\
&\ \text{curl}\left(\text{curl}\ b\right) \in L^2(\Omega) \bigg \rbrace. \nonumber
\end{align}

While the more general space $W^k(\Omega)$ in (\ref{937}), containing also functions which are not divergence-free, serves as a test function space for the induction equation at the discrete time $k\Delta t$, the space $Y^k(\Omega)$ will be the space in which we construct the magnetic induction at time $k \Delta t$. As for functions $b \in Y^k(\Omega)$ it holds $\text{curl}(\text{curl}b)=\Delta b$, both of these spaces can be equipped with the $H^{2}$-norm.

\begin{proposition}
\label{Existence}

Let all the assumptions of Theorem \ref{mainresult} be satisfied and $\Delta t > 0$. Let further $g_{\Delta t}^k$ and $J_{\Delta t}^k$ be given by (\ref{840}) for any $k=0,...,\frac{T}{\Delta t}$ and assume in addition that
\begin{align}
\rho _0 \in H^{1}(\Omega), \quad u_0,\ B_0 \in H^{2}(\Omega). \nonumber
\end{align}

Then, for all $k=1,...,\frac{T}{\Delta t}$, there exist functions $\chi_{\Delta t,k} \in C([(k-1)\Delta t,k \Delta t];L^p_\text{loc}(\mathbb{R}^3))$, $1 \leq p < \infty$ and
\begin{align}
\rho_{\Delta t}^k \in H^{1}(\Omega),\quad \underline{\rho} \leq \rho_{\Delta t}^k \leq \overline{\rho},\quad \quad u_{\Delta t}^k \in V_0^2(\Omega),\quad B_{\Delta t}^k \in Y^k(\Omega) \label{852}
\end{align}

which satisfy the variational equations (\ref{236})--(\ref{116}) for all test functions $\Theta \in \mathcal{D}([(k-1)\Delta t,k\Delta t]\times \mathbb{R}^3)$, $\psi \in H^{1}(\Omega)$, $\phi \in V_0^2(\Omega)$ and $b \in W^k(\Omega)$.

\end{proposition}

\textbf{Proof}

\bigskip

We consider some discrete time index $k \in \left\lbrace 1,...,\frac{T}{\Delta t} \right\rbrace$ and assume that the proposition is already proved for all $l=1,...,k-1$.

\bigskip

\textbf{Step 1:} The existence of a solution $\rho_{\Delta t}^k \in H^{1}(\Omega)$ to (\ref{8}) follows immediately from the Lax-Milgram Lemma. In order to show the upper bound of $\rho_{\Delta t}^k$ in (\ref{852}), we follow \cite[Section 7.6.5]{novotnystraskraba} and set $r_{\Delta t}^k := \rho_{\Delta t}^k - \overline{\rho}$. Since $\rho_{\Delta t}^k$ is a solution to (\ref{8}) and $\overline{\rho}$ is a constant it follows that
\begin{align}
-\int _\Omega \frac{r_{\Delta t}^k - r_{\Delta t}^{k-1}}{\Delta t}\psi  dx = \int_\Omega u_{\Delta t}^{k-1} \cdot \nabla r _{\Delta t}^k \psi + \epsilon \nabla r^k_{\Delta t} \cdot \nabla \psi \ dx \nonumber
\end{align}

for all $\psi \in H^1(\Omega)$. Testing this equation by $\max \lbrace r_{\Delta t}^k,0 \rbrace \in H^1(\Omega)$ and using that $\operatorname{div}u_{\Delta t}^{k-1} = 0$ we infer that
\begin{align}
\int _\Omega \frac{\left|\max\left\lbrace r_{\Delta t}^k, 0 \right\rbrace\right|^2}{\Delta t} + \epsilon \left| \nabla \max \lbrace r_{\Delta t}^k,0\rbrace \right|^2\ dx = \int _\Omega \frac{r_{\Delta t}^{k-1} \max\left\lbrace r_{\Delta t}^k, 0 \right\rbrace}{\Delta t}\ dx. \nonumber
\end{align}

Since $r_{\Delta t}^{k-1} = \rho_{\Delta t}^{k-1} - \overline{\rho} \leq 0$, the right-hand side of this equation is nonpositive and we infer that $r_{\Delta t}^k \leq 0$, i.e.\ $\rho_{\Delta t}^k \leq \overline{\rho}$. Arguing similarly for the lower bound, we arrive at the estimates for $\rho_{\Delta t}^k$ in (\ref{852}).

\bigskip

\textbf{Step 2:} As in \cite{gigli}, we consider the initial value problem
\begin{align}
\frac{\partial X_{\Delta t}^{\Pi _{\Delta t}^{k-1}}(s;t,x)}{\partial t} = \Pi _{\Delta t}^{k-1} \left( X_{\Delta t}^{\Pi _{\Delta t}^{k-1}}(s;t,x) \right), \quad \quad X_{\Delta t}^{\Pi _{\Delta t}^{k-1}}(s;s,x) = x,\quad x \in \mathbb{R}^3,\quad s,t \in \mathbb{R}, \label{143}
\end{align}

where $t$ represents the time variable and $s$ the initial time. Since $\Pi _{\Delta t}^{k-1}$ is constant in time and a rigid velocity field by (\ref{222}), it is in particular Lipschitz-continuous. Then by the theory of ordinary differential equations, (\ref{143}) defines a unique mapping
\begin{equation}
X_{\Delta t}^{\Pi _{\Delta t}^{k-1}}: \mathbb{R} \times \mathbb{R} \times \mathbb{R}^3 \rightarrow \mathbb{R}^3. \label{144}
\end{equation}

We set
\begin{equation}
\chi_{\Delta t,k}(t,x) := \chi_{\Delta t, k-1} \left((k-1)\Delta t, X_{\Delta t}^{\Pi_{\Delta t}^{k-1}} \left(t;(k-1)\Delta t,x\right) \right)\quad \text{for } t \in [(k-1)\Delta t,k\Delta t] \label{944}
\end{equation}

and infer from \cite[Theorem III.2]{dipernalions} that this is the (unique renormalized) solution to the transport equation (\ref{236}).

\bigskip

\textbf{Step 3:} We consider the operator
\begin{align}
A: V_0^2(\Omega) \rightarrow \left(V_0^2(\Omega)\right)^*,\quad \left\langle Au,v \right\rangle _{\left(V_0^2(\Omega)\right)^*\times V_0^2(\Omega)} :=& \int _\Omega \left( \frac{\rho _{\Delta t}^k u}{\Delta t}\right) \cdot v + \text{div}\left( \rho _{\Delta t}^{k} u_{\Delta t}^{k-1} \otimes u \right)\cdot v \nonumber \\
+& 2\nu \mathbb{D}(u) : \nabla v + \epsilon \left( \nabla u \nabla \rho_{\Delta t}^k \right) \cdot v + \epsilon (\Delta u) \cdot (\Delta v)\ dx \nonumber
\end{align}

for $u,v \in V_0^2(\Omega)$. Because of the regularization term, $A$ is coercive on $V_0^2(\Omega)$. Further, the bilinear form $\left\langle A\cdot,\cdot \right\rangle _{\left(V_0^2(\Omega)\right)^*,V_0^2(\Omega)}$ is bounded on $V_0^2(\Omega)$ and hence the Lax-Milgram Lemma again implies the existence of $u_{\Delta t}^k \in V_0^2(\Omega)$ satisfying (\ref{57}).

\bigskip

\textbf{Step 4:} We introduce
\begin{align}
\tilde{A}:Y^k(\Omega)& \rightarrow \left( Y^k(\Omega) \right)^*, \nonumber \\
\left\langle \tilde{A}(B),b \right\rangle _{\left(Y^k(\Omega)\right)^* \times Y^k(\Omega)} :=& \int _\Omega  \frac{B}{\Delta t} \cdot b + \epsilon \text{curl}\left(\text{curl} B \right) \cdot \text{curl}\left(\text{curl} b\right) \nonumber \\
&+ \left[ \frac{1}{\sigma \mu} \text{curl}B + \frac{\epsilon}{\mu ^2} \left| \text{curl} B \right|^2\text{curl} B \right] \cdot \text{curl} b\ dx. \nonumber
\end{align}

Clearly, $\tilde{A}$ is coercive:
\begin{align}
\left\langle \tilde{A}(B),B \right\rangle_{\left(Y^k\right)^* \times Y^k} &\geq \frac{1}{\Delta t} \| B \|_{L^2(\Omega)}^2 + \epsilon\| \Delta B \|_{L^2(\Omega)}^2 \geq c\|B\|^2_{H^{2}(\Omega)} = c \|B\|^2_{Y^k}. \nonumber
\end{align}

Further, if $B_n \rightharpoonup B$ in $Y^k(\Omega)$, then from the Rellich-Kondrachov embedding we know $B_n \rightarrow B$ in $W^{1,4}(\Omega)$, which again gives us weak continuity of $\tilde{A}$. Coercivity and weak continuity imply surjectivity of $\tilde{A}$ (see for example \cite[Theorem 1.2]{francu}) and so we infer the existence of a solution $B_{\Delta t}^k \in Y^k(\Omega)$ to (\ref{116}) for all $b \in Y^k(\Omega)$. Moreover, by the Helmholtz-decomposition \cite[Theorem 4.2]{maxwellinequality}, we find for any $b \in W^k(\Omega)$ some functions
\begin{align}
q \in H^{1}(\Omega),\quad \quad w \in\lbrace \tilde{w} \in L^2(\Omega):\ \text{curl} \tilde{w} \in L^2(\Omega),\quad \nabla \cdot \tilde{w} = 0\quad \text{in } \Omega,\quad \tilde{w} \times \operatorname{n}|_{\partial \Omega} = 0 \rbrace \nonumber
\end{align}

such that
\begin{align}
b = \nabla q + \operatorname{curl}w. \nonumber
\end{align}

It is easy to see that $\operatorname{curl}w \in Y^k(\Omega)$, i.e. $\operatorname{curl}w$ is an admissible test function in (\ref{116}). From the fact that $\operatorname{curl} (\nabla q) = 0$ and $\nabla \cdot B_{\Delta t}^k = \nabla \cdot B_{\Delta t}^{k-1} = 0$ we conclude that (\ref{116}) even holds true for all $b \in W^k(\Omega)$.

$\hfill \Box$

\begin{remark}
\label{isometries}
For any fixed $s,t \in \mathbb{R}$ the mapping (\ref{144}) is an isometry. Indeed, from $\Pi _{\Delta t}^{k-1}$ being a rigid velocity field and the ordinary differential equation (\ref{143}), it follows that
\begin{equation}
\frac{\partial}{\partial t} \left| X_{\Delta t}^{\Pi _{\Delta t}^{k-1}}(s;t,x) - X_{\Delta t}^{\Pi _{\Delta t}^{k-1}}(s;t,y) \right|^2 = 0 \label{283}
\end{equation}

for any $x,y \in \mathbb{R}^3$. 
\end{remark}

In the remainder of this section we derive an energy inequality for our discrete solution. To this end we extend, without loss of generality,
\begin{align}
u_{\Delta t}^l(x) = 0,\quad \rho_{\Delta t}^l(x) = \underline{\rho},\quad \quad \forall x \in \mathbb{R}^3 \setminus \Omega,\ l = 0,...,\frac{T}{\Delta t}. \label{834}
\end{align}

We fix some $k \in \left\lbrace 1,...,\frac{T}{\Delta t} \right\rbrace$. For arbitrary $l \leq k$ we test the continuity equation (\ref{8}) at the discrete time $l \Delta t$ by $\frac{1}{2}|u_{\Delta t}^l|^2$ and subtract the result from the momentum equation (\ref{57}), also at time $l \Delta t$, tested by $u_{\Delta t}^l$. This yields
\begin{eqnarray}
&& \int _\Omega \frac{1}{2\Delta t} \rho _{\Delta t}^l |u_{\Delta t}^l|^2 - \frac{1}{2\Delta t} \rho _{\Delta t}^{l-1} |u_{\Delta t}^{l-1}|^2 + 2\nu |\nabla u_{\Delta t}^l|^2 + \frac{1}{\eta} \rho _{\Delta t}^{l-1}\chi_{\Delta t}^{l}(u_{\Delta t}^{l-1} - \Pi _{\Delta t}^{l-1} ) \cdot u_{\Delta t}^l + \epsilon |\Delta u_{\Delta t}^l|^2\ dx \nonumber \\
&\leq& \int _\Omega \rho _{\Delta t}^{l-1} g_{\Delta t}^l \cdot u_{\Delta t}^l + \frac{1}{\mu}(\text{curl}B_{\Delta t}^{l-1} \times B_{\Delta t}^{l-1}) \cdot u_{\Delta t}^l\ dx. \label{17}
\end{eqnarray}

Next, we test the magnetic equation (\ref{116}) at time $l \Delta t$ by $\frac{1}{\mu}B_{\Delta t}^l$ and estimate
\begin{align}
&\frac{1}{2 \mu \Delta t} \left[ \|B_{\Delta t}^l\|_{L^2(\Omega)}^2 - \|B_{\Delta t}^{l-1}\|_{L^2(\Omega)}^2 \right] + \frac{1}{\mu}\int _\Omega \frac{1}{\sigma \mu} |\text{curl}B_{\Delta t}^l|^2\ dx \nonumber \\
\leq& - \int _\Omega \frac{\epsilon}{\mu ^3} |\text{curl} B_{\Delta t}^l|^4 + \frac{\epsilon}{\mu} \left| \Delta B_{\Delta t}^l \right|^2 - \frac{1}{\mu} (u_{\Delta t}^l \times B_{\Delta t}^{l-1})\cdot \text{curl}B_{\Delta t}^l - \frac{1}{\sigma \mu} J_{\Delta t}^l \cdot \text{curl}B_{\Delta t}^l\ dx. \label{340}
\end{align}

Adding this to (\ref{17}) and summing over all $l\leq k$, we infer
\begin{align}
&\frac{1}{2\Delta t}\underline{\rho}\|u_{\Delta t}^k\|_{L^2(\Omega)}^2 - \frac{1}{2\Delta t}\overline{\rho}\|u_{\Delta t}^0\|_{L^2(\Omega)}^2 + \sum_{l=1}^k\left( 2\nu \| \nabla u_{\Delta t}^l \|_{L^2(\Omega)}^2 + \epsilon \| \Delta u_{\Delta t}^l \|_{L^2(\Omega)}^2 \right) \nonumber \\
+& \frac{1}{2 \mu \Delta t} \left( \| B_{\Delta t}^k \|_{L^2(\Omega)}^2 - \| B_{\Delta t}^0 \|_{L^2(\Omega)}^2 \right) + \sum _{l=1}^k \left( \frac{1}{\sigma \mu} \| \text{curl}B_{\Delta t}^l \|_{L^2(\Omega)}^2 + \frac{\epsilon }{\mu ^3} \| \text{curl} B_{\Delta t}^l \|_{L^4(\Omega)}^4 + \frac{\epsilon}{\mu} \left\| \Delta B_{\Delta t}^l \right\|_{L^2(\Omega)}^2 \right) \nonumber \\
\leq& \sum _{l=1}^k \int _\Omega -  \frac{1}{\eta} \rho _{\Delta t}^{l-1}\chi_{\Delta t}^{l}(u_{\Delta t}^{l-1} - \Pi _{\Delta t}^{l-1} ) \cdot u_{\Delta t}^l + \rho _{\Delta t}^{l-1}g_{\Delta t}^l \cdot u_{\Delta t}^l + \frac{1}{\mu}(\text{curl}B_{\Delta t}^{l-1} \times B_{\Delta t}^{l-1}) \cdot u_{\Delta t}^l \nonumber \\
+& \frac{1}{\mu} (u_{\Delta t}^l \times B_{\Delta t}^{l-1})\cdot \text{curl}B_{\Delta t}^l + \frac{1}{\sigma \mu} J_{\Delta t}^l \cdot \text{curl}B_{\Delta t}^l dx. \label{18}
\end{align}

In order to estimate the right-hand side here, we need
\begin{equation}
\| \Pi_{\Delta t}^{l-1} \|_{L^2(\Omega)} \leq c \| u_{\Delta t}^{l-1} \|_{L^2(\Omega)},\label{228}
\end{equation}

which can be proved in the following way: We distinguish between two cases, the first one being $\text{supp}\chi_{\Delta t}^{l-1} \bigcap \Omega = \emptyset$. Then, as $u_{\Delta t}^{l-1}=0$ outside of $\Omega$, both sides of (\ref{228}) are equal to zero and so the inequality is trivially satisfied. For the second case, $\text{supp}\chi_{\Delta t}^{l-1} \bigcap \Omega \neq \emptyset$, we note that from (\ref{144}) being an isometry and by (\ref{944}) it follows the existence of a compact set $K$ independent of $l$ and $\Delta t$ such that in this case $\text{supp}\chi_{\Delta t}^{l-1} \subset K$. This allows us to reduce the integrals over $\mathbb{R}^3$ in (\ref{835}) - (\ref{837}) to integrals over $K$. Moreover, as $\rho _{\Delta t}^l \geq \underline{\rho}$ was extended by $\underline{\rho}$ on $\mathbb{R}^3 \setminus \Omega$, we know
\begin{equation}
\int _{\mathbb{R}^3} \rho _{\Delta t}^l \chi_{\Delta t}^{l-1}\ dx \geq \underline{\rho} |S|>0. \label{216}
\end{equation}

Thus we can estimate
\begin{align}
\left| a_{\Delta t}^{l-1} \right| = \left| \frac{\int _{\mathbb{R}^3}\rho _{\Delta t}^l \chi_{\Delta t}^{l-1}x\ dx}{\int _{\mathbb{R}^3} \rho _{\Delta t}^l \chi_{\Delta t}^{l-1}l\ dx } \right| \leq c \left| \int _K \rho _{\Delta t}^lx\ dx \right| \leq c \label{935}
\end{align}

with $c$ independent of $l$ and $\Delta t$. By similar computations, c.f.\ also \cite[Section 3.2]{cottetmaitre} and the proof of \cite[Lemma 4]{alternativerigidbodies}, we obtain
\begin{align}
\left| \left( u_G \right)_{\Delta t}^{l-1} \right| \leq c \left\| u_{\Delta t}^{l-1} \right\|_{L^2(\Omega)},\quad \left| \omega_{\Delta t}^{l-1} \right| \leq c \left\| u_{\Delta t}^{l-1} \right\|_{L^2(\Omega)}, \quad v \cdot \left( I_{\Delta t}^{l-1}v\right) \geq c|u|^2 \quad \forall v \in \mathbb{R}^3, \label{849}
\end{align}

where the last inequality uses that since (\ref{144}) is an isometry, one can find for any $\Delta t > 0$, $l = 1,...,\frac{T}{\Delta t}$ some ball $B_r(l,\Delta t) \subset \mathbb{R}^3$ with radius $r>0$ independent of $l$ and $\Delta t$ such that $B_r(l,\Delta t) \subset S(\chi_{\Delta t}^{l-1})$. Thus, (\ref{228}) is also satisfied in the second case. Now, exploiting (\ref{228}) and applying Young's inequality, the right-hand side of (\ref{18}) can be bounded by
\begin{align}
& \sum_{l=1}^k \left[ \frac{\overline{\rho}}{2\eta} \|u_{\Delta t}^{l-1}\|_{L^2(\Omega)}^2 + \frac{\overline{\rho}}{2\eta} \|u_{\Delta t}^l\|_{L^2(\Omega)}^2 + \frac{c^2 \overline{\rho}}{2} \|g\|_{L^\infty (\Omega)}^2 + \frac{\overline{\rho}}{2} \| u_{\Delta t}^l \|_{L^2(\Omega)}^2\right.\nonumber \\
& \left.+ \frac{1}{\mu}\| \text{curl} B_{\Delta t}^{l-1} \|_{L^4(\Omega)}\|B_{\Delta t}^{l-1}\|_{L^4 (\Omega)}\| u_{\Delta t}^l \|_{L^2(\Omega)} + \frac{c^2}{2\sigma \mu} \| J_{\Delta t}^l \|_{L^\infty (\Omega)}^2 + \frac{1}{2\sigma \mu} \|\text{curl}B_{\Delta t}^l\|_{L^2(\Omega)}^2 \right. \nonumber \\
& \left. + \frac{1}{\mu} \| \text{curl} B_{\Delta t}^{l} \|_{L^4(\Omega)}\|B_{\Delta t}^{l-1}\|_{L^4 (\Omega)}\| u_{\Delta t}^l \|_{L^2(\Omega)} \right]. \label{838}
\end{align}

Using the Poincaré-type estimate
\begin{equation}
\left\| B_{\Delta t}^{l-1} \right\|_{L^4(\Omega)} \leq c\left\| \text{curl}B_{\Delta t}^{l-1} \right\|_{L^4(\Omega)}, \nonumber
\end{equation}

c.f.\ \cite[Corollary 3.4]{friedrichsinequality}, we further estimate
\begin{eqnarray}
\frac{1}{\mu}\| \text{curl} B_{\Delta t}^{l} \|_{L^4(\Omega)}\|B_{\Delta t}^{l-1}\|_{L^4(\Omega)}\| u_{\Delta t}^l \|_{L^2(\Omega)} &\leq& \frac{c^2\mu}{\epsilon} \| u_{\Delta t}^l \|_{L^2(\Omega)}^2 + \frac{\epsilon}{8\mu ^3} \| \text{curl} B_{\Delta t}^{l-1} \|_{L^4(\Omega)}^4 + \frac{\epsilon}{8 \mu ^3} \| \text{curl} B_{\Delta t}^l \|_{L^4(\Omega)}^4, \nonumber \\
\frac{1}{\mu}\| \text{curl} B_{\Delta t}^{l-1} \|_{L^4(\Omega)}\|B_{\Delta t}^{l-1}\|_{L^4(\Omega)}\| u_{\Delta t}^l \|_{L^2(\Omega)} &\leq& \frac{c^2\mu}{\epsilon} \| u_{\Delta t}^l \|_{L^2(\Omega)}^2 + \frac{\epsilon}{4\mu^3}\| \text{curl} B_{\Delta t}^{l-1} \|_{L^4(\Omega)}^4. \nonumber
\end{eqnarray}

Consequently, we can absorb several quantities from (\ref{838}), including - for $\Delta t$ sufficiently small - the $L^2$-norm of $u_{\Delta t}^k$, into the left-hand side of (\ref{18}) and obtain
\begin{eqnarray}
&&\frac{\underline{\rho}}{4\Delta t}\left\|u_{\Delta t}^k\right\|_{L^2(\Omega)}^2 + \sum_{l=1}^k\left( 2\nu \left\| \nabla u_{\Delta t}^l \right\|_{L^2(\Omega)}^2 + \epsilon \left\| \Delta u_{\Delta t}^l \right\|^2_{L^2(\Omega)} \right) + \frac{1}{2\mu \Delta t} \left\| B_{\Delta t}^k \right\|_{L^2(\Omega)}^2 \nonumber \\
&&+ \sum _{l=1}^k \left( \frac{1}{2\sigma \mu} \left\| \text{curl} B_{\Delta t}^l \right\|^2_{L^2(\Omega)} + \frac{\epsilon }{2\mu ^3} \left\| \text{curl} B_{\Delta t}^l \right\|_{L^4(\Omega)}^4 + \frac{\epsilon}{\mu} \left\| \Delta B_{\Delta t}^l \right\|_{L^2(\Omega)}^2 \right) \nonumber \\
&\leq& \frac{\overline{\rho}}{2\Delta t}\left\|u_{\Delta t}^0\right\|_{L^2(\Omega)}^2 + \frac{1}{2\mu \Delta t}\left\|B_{\Delta t}^0\right\|_{L^2(\Omega)}^2 + \frac{T c^2 \overline{\rho}}{2\Delta t} \left\|g\right\|_{L^\infty ((0,T)\times \Omega)}^2 + \frac{T c^2}{2\sigma \mu \Delta t} \left\|\text{curl}J\right\|_{L^\infty ((0,T)\times \Omega)}^2 \nonumber \\
&&+ \frac{\epsilon}{8\mu ^3}\left\| \text{curl} B_{\Delta t}^0 \right\|_{L^4(\Omega)}^4 + \frac{\epsilon}{4\mu ^3}\left\| \text{curl} B_{\Delta t}^0 \right\|_{L^4(\Omega)}^4 + \frac{\overline{\rho}}{2\eta} \left\| u_{\Delta t}^0 \right\|_{L^2(\Omega)}^2 + \sum _{l=1}^{k-1} \left[ \frac{\overline{\rho}}{\eta} \left\| u_{\Delta t}^l \right\|_{L^2(\Omega)}^2 + \frac{\overline{\rho}}{2} \left\| u_{\Delta t}^l \right\|_{L^2(\Omega)}^2 \right. \nonumber \\
&&+ \left.\frac{2c^2\mu}{\epsilon} \left\| u_{\Delta t}^l \right\|_{L^2(\Omega)}^2 \right]. \label{19}
\end{eqnarray}

Hence, from the discrete Gronwall estimate (c.f.\ \cite[(1.67)]{roubicek}), we infer the bound
\begin{eqnarray}
&& \left\|u_{\Delta t}^k\right\|_{L^2(\Omega)}^2 + \Delta t\sum_{l=1}^k\left( \left\| \nabla u_{\Delta t}^l \right\|_{L^2(\Omega)}^2 + \left\| \Delta u_{\Delta t}^l \right\|_{L^2(\Omega)}^2 \right) + \left\| B_{\Delta t}^k \right\|_{L^2(\Omega)}^2\nonumber \\
&&+ \Delta t \sum _{l=1}^k \left( \left\| \text{curl}{B}_{\Delta t}^l \right\|^2_{L^2(\Omega)} + \left\| \text{curl} B_{\Delta t}^l \right\|_{L^4(\Omega)}^4 + \left\| \Delta B_{\Delta t}^l \right\|_{L^2(\Omega)}^2 \right) \nonumber \\
&\leq& c(u^0, B^0, \overline{\rho}, \underline{\rho}, c, g, J, \sigma, \mu, \nu, \epsilon, \eta , T) \quad \forall k = 1,...,\frac{T}{\Delta t} \label{20}
\end{eqnarray}

uniformly in $\Delta t$ and $k$.

{\centering \section{Limit passage with respect to \texorpdfstring{$\Delta t \rightarrow 0$}{}} \label{rothelimit} \par 
}
We now want to pass to the limit in the time discretization, i.e.\ $\Delta t \to 0$. To do so, we introduce piecewise constant as well as piecewise affine interpolants of our functions defined so far only in the discrete time points. Namely, for the time-independent quantities $f_{\Delta t}^k$ defined for $k=0,...,\frac{T}{\Delta t}$ we set
\begin{align}
f_{\Delta t}(t) &:= \left( \frac{t}{\Delta t} - (k-1)\right)f_{\Delta t}^k + \left(k - \frac{t}{\Delta t}\right) f_{\Delta t}^{k-1}\ \ \ &\text{for }& (k-1)\Delta t < t \leq k \Delta t,\quad k=1,...,\frac{T}{\Delta t}, \label{21} \\
\overline{f}_{\Delta t}(t) &:= f_{\Delta t}^k \ \ \ &\text{for }& (k-1)\Delta t < t \leq k \Delta t,\quad k=0,...,\frac{T}{\Delta t}, \label{22} \\
\overline{f}'_{\Delta t}(t) &:= f_{\Delta t}^{k-1} \ \ \ &\text{for }& (k-1)\Delta t < t \leq k \Delta t,\quad k=1,...,\frac{T}{\Delta t}. \label{23}
\end{align}

We will use the same notation also for the interpolation of the discrete momentum function $(\rho u)_{\Delta t}^k := \rho _{\Delta t}^k u_{\Delta t}^k$, $k = 0,...,\frac{T}{\Delta t}$. Regarding the solution to the transport equation on $[0,T]$, we glue together the already time-dependent functions $\chi_{\Delta t,k}$, defined on the intervals $[(k-1)\Delta t, k \Delta t]$. More specifically, we set
\begin{align}
\chi_{\Delta t}(t) := \chi_{\Delta t,k}(t) \quad \quad \text{for } (k-1)\Delta t < t \leq k\Delta t,\quad k = 1,...,\frac{T}{\Delta t}. \nonumber
\end{align}

By the construction of $\chi_{\Delta t, k}$ in Proposition \ref{Existence} it holds $\chi_{\Delta t} \in C([0,T];L^p_{\text{loc}}(\mathbb{R}^3))$, $1 \leq p < \infty$, and $\chi_{\Delta t}$ is the solution to
\begin{equation}
-\int_0^{T} \int_{\mathbb{R}^3 } \chi_{\Delta t} \partial_t \Theta dxdt - \int _{\mathbb{R}^3 } \chi _0 \Theta (0,x)\ dx = \int_0^T \int_{\mathbb{R}^3} \left( \chi_{\Delta t} \overline{\Pi}'_{\Delta t} \right) \cdot \nabla \Theta \ dxdt \label{926}
\end{equation}

for any $\Theta \in \mathcal{D}([0,T)\times \mathbb{R}^3)$. According to the transport theory by DiPerna and Lions, c.f.\ \cite[Thoerem III.2]{dipernalions}, this solution is unique and can be represented by
\begin{equation}
\chi_{\Delta t}(t,x) := \chi_0 \left( X_{\Delta t}^{\overline{\Pi}'_{\Delta t}} (t;0,x) \right)\ \ \ \text{for } t \in [0,T]. \label{119}
\end{equation}
Here $X_{\Delta t}^{\overline{\Pi}'_{\Delta t}}$ denotes the unique solution to the initial value problem
\begin{align}
\frac{\partial X_{\Delta t}^{\overline{\Pi}'_{\Delta t}} (s;t,x)}{\partial t} = \overline{\Pi}'_{\Delta t} \left(t, X_{\Delta t}^{\overline{\Pi}'_{\Delta t}} (s;t,x) \right), \quad \quad X_{\Delta t}^{\overline{\Pi}'_{\Delta t}} (s;s,x) = x,\quad x \in \mathbb{R}^3,\quad s,t \in [0,T], \label{234}
\end{align}
given by the Carathéodory theorem \cite[Theorem 1.45]{roubicek}. By the uniqueness of this solution, the function $X_{\Delta t}^{\overline{\Pi}'_{\Delta t}}$ can also be written as a composition of the mappings (\ref{144}). In particular, by the corresponding property of those functions (c.f.\ Remark \ref{isometries}), the mapping
\begin{equation}
x \rightarrow X_{\Delta t}^{\overline{\Pi}_{\Delta t}'}(s;t,x),\ \ \ s,t \in [0,k\Delta t]\label{284}
\end{equation}

is an isometry from $\mathbb{R}^3$ to $\mathbb{R}^3$. The a-priori estimate (\ref{20}) translates to the following uniform bounds for the above defined interpolants
\begin{align}
\left\| u_{\Delta t} \right\|_{L^\infty (0,T;L^2(\Omega))} + \left\| \overline{u}_{\Delta t} \right\|_{L^\infty (0,T;L^2(\Omega))} + \left\| \overline{u}'_{\Delta t} \right\|_{L^\infty (0,T;L^2(\Omega))} &\leq c, \label{29} \\
\left\| u_{\Delta t} \right\|_{L^2 (0,T;H^{2}(\Omega))} + \left\| \overline{u}_{\Delta t} \right\|_{L^2 (0,T;H^{2}(\Omega))} + \left\| \overline{u}'_{\Delta t} \right\|_{L^2 (0,T;H^{2}(\Omega))} &\leq c, \label{543} \\
\left\| B_{\Delta t} \right\|_{L^\infty (0,T;L^2(\Omega))} + \left\| \overline{B}_{\Delta t} \right\|_{L^\infty (0,T;L^2(\Omega))} + \left\| \overline{B}'_{\Delta t} \right\|_{L^\infty (0,T;L^2(\Omega))} &\leq c, \label{31} \\
\left\| B_{\Delta t} \right\|_{L^2 (0,T;H^{2}(\Omega))} + \left\| \overline{B}_{\Delta t} \right\|_{L^2 (0,T;H^{2}(\Omega))} + \left\| \overline{B}'_{\Delta t} \right\|_{L^2(0,T;H^{2}(\Omega))} &\leq c, \label{344} \\
\left\| \operatorname{curl}B_{\Delta t} \right\|_{L^4 ((0,T)\times \Omega)} + \left\| \operatorname{curl}\overline{B}_{\Delta t} \right\|_{L^4 ((0,T)\times \Omega)} + \left\| \operatorname{curl}\overline{B}'_{\Delta t} \right\|_{L^4 ((0,T)\times \Omega)} &\leq c. \label{920}
\end{align}

These bounds allow us to find functions
\begin{align}
B &\in \bigg\lbrace b \in L^\infty \left(0,T;L^2(\Omega)\right) \bigcap L^2\left(0,T;V^2(\Omega)\right):\ b \cdot \text{n}|_{\partial \Omega} = 0 \bigg\rbrace \label{825} \\
u &\in  L^\infty \left(0,T;L^2(\Omega)\right) \bigcap L^2\left(0,T;V_0^2(\Omega)\right) \label{880}
\end{align}

such that for selected non-relabeled subsequences
\begin{align}
\overline{B}'_{\Delta t},\ \overline{B}_{\Delta t},\ B_{\Delta t} \buildrel\ast\over\rightharpoonup B \ \ \ &\text{in } L^\infty \left(0,T;L^2\left(\Omega \right)\right),\ \ \ &\overline{B}'_{\Delta t},\ \overline{B}_{\Delta t},\ B_{\Delta t}& \rightharpoonup B \ \ \ \text{in } L^2\left(0,T;H^{2}\left(\Omega \right)\right), \label{38} \\
\overline{u}'_{\Delta t},\ \overline{u}_{\Delta t},\ u_{\Delta t} \buildrel\ast\over\rightharpoonup u \ \ \ &\text{in } L^\infty \left(0,T;L^2\left(\Omega \right)\right),\ \ \ &\overline{u}'_{\Delta t},\ \overline{u}_{\Delta t},\ u_{\Delta t}& \rightharpoonup u \ \ \ \text{in } L^2\left(0,T;H^{2}\left(\Omega \right)\right). \label{796}
\end{align}

The equality between the weak limits of $\overline{B}'_{\Delta t},\ \overline{B}_{\Delta t},\ B_{\Delta t}$ and $\overline{u}'_{\Delta t},\ \overline{u}_{\Delta t},\ u_{\Delta t}$ is given by Lemma \ref{equalityofrothelimits}. The inclusions (\ref{825}) and (\ref{880}) follow from the fact that $\overline{u}_{\Delta t}(t) \in V_0^2$ for any $t \in [0,T]$ and $\overline{B}_{\Delta t}(t) \in Y^k$ for any $t \in ((k-1)\Delta t,k\Delta t]$. Moreover, for the discretized external forces $\overline{J}_{\Delta t}$ and $\overline{g}_{\Delta t}$ it follows
\begin{align}
\overline{J}_{\Delta t} \rightarrow J \quad \text{in }L^p((0,T)\times \Omega),\quad \quad \overline{g}_{\Delta t} \rightarrow g \quad \text{in }L^p((0,T)\times \Omega) \quad \quad \quad \forall 1 \leq p < \infty, \nonumber
\end{align}

directly from their definition in (\ref{840}), c.f.\ \cite[Lemma 8.7]{roubicek}. \\

\subsection{Characteristic function} \label{rothelimitcharacteristics}

The fact that it still holds $\overline{\Pi}_{\Delta t}^{l-1}=0$ whenever $\text{supp}\chi_{\Delta t}^l \bigcap \Omega = \emptyset$ and (\ref{935}), (\ref{849}), (\ref{20}) in the other case imply the condition (\ref{841}) from Lemma \ref{strongconvergencecharacteristics} and in particular we get the existence of a function $\Pi \in L^\infty(0,T;W^{1,\infty}_\text{loc} (\mathbb{R}^3))$ such that
\begin{equation}
\overline{\Pi}'_{\Delta t} \buildrel\ast\over\rightharpoonup \Pi \ \ \ \text{in } L^\infty \left(0,T;W^{1,\infty}_\text{loc} (\mathbb{R}^3)\right),\quad \quad \Pi(t,x) = v(t) + w(t) \times x,\quad \quad v,w \in L^\infty(0,T). \label{153}
\end{equation}

In (\ref{221}) we will characterize the limit function $\Pi$ more specifically through the density, the velocity and the characteristic function, similar to (\ref{222}). By the relations (\ref{926}) and (\ref{234}) we also have the conditions (\ref{845}) and (\ref{850}) of Lemma \ref{strongconvergencecharacteristics} which, in combination with Remark \ref{strongconvergencecharacteristicsremark}, implies that
\begin{align}
\chi_{\Delta t} \rightarrow \chi \quad &\text{in } C\left([0,T];L^p\left(\mathbb{R}^3\right)\right)\quad \forall 1 \leq p < \infty, \quad \quad \chi (t,x) = \chi_0\left(X^{\Pi}(t;0,x)\right), \label{128} \\
X_{\Delta t}^{\overline{\Pi}'_{\Delta t}} \rightarrow X^{\Pi} \ \ \ &\text{in } C\left([0,T]\times [0,T];C_\text{loc}\left(\mathbb{R}^3 \right)\right), \label{117}
\end{align}

where $\chi$ and $X^\Pi$ are the unique solutions of
\begin{align}
-\int_0^T \int_{\mathbb{R}^3 } \chi \partial_t \Theta dxdt - \int _{\mathbb{R}^3 } \chi _0 \Theta (0,x)\ dx &= \int_0^T \int_{\mathbb{R}^3} \left(\chi \Pi \right) \cdot \nabla \Theta \ dxdt \quad \forall \Theta \in \mathcal{D} \left([0,T)\times \mathbb{R}^3\right), \label{288} \\
\frac{dX^{\Pi}(s;t,x)}{dt} &= \Pi \left(t, X^{\Pi} (s;t,x) \right), \quad \quad X^{\Pi}(s;s,x)= x \label{864}
\end{align}

respectively. From (\ref{128}) it also follows that
\begin{align}
\overline{\chi}_{\Delta t},\ \overline{\chi}'_{\Delta t} \rightarrow \chi \quad &\text{in } C\left([0,T];L^p\left(\mathbb{R}^3\right)\right)\quad \forall 1 \leq p < \infty, \label{919}
\end{align}

which is obtained in the same way as the similar statement in \cite[Lemma 8.7]{roubicek}.

\subsection{Induction equation}\label{inductionrothelimit}

In the passage to the limit of the induction equation (\ref{265}) below, we consider test functions from the space $Y(\chi,T)$ which are curl-free in a neighbourhood of the solid region in the limit. To see that this is possible, let us choose an arbitrary $\gamma > 0$ and denote by $S_{\gamma}(\chi (t))$ and $S^{\gamma}(\chi (t))$ the $\gamma$-neighbourhood and the “$\gamma$-kernel” of $S(\chi(t))$ respectively, i.e.
\begin{align}
S^{\gamma}(\chi (t)):= \left\lbrace x \in \mathbb{R}^3:\ \operatorname{dist}\left(x,S(\chi (t)) \right) < \gamma \right\rbrace, \quad S_{\gamma}(\chi (t)):= \left\lbrace x \in S(\chi (t)):\ \operatorname{dist}\left(x,\partial S(\chi (t)) \right) > \gamma \right\rbrace. \nonumber
\end{align}

From the uniform convergence (\ref{117}) of $X_{\Delta t}^{\overline{\Pi}'_{\Delta t}}$ and the relation (\ref{119}) between $X_{\Delta t}^{\overline{\Pi}'_{\Delta t}}$ and the characteristic functions $\chi_{\Delta t}$, it follows the existence of some $\delta (\gamma) > 0$ such that
\begin{align}
S_\gamma (\chi (t)) \subset S \left( \overline{\chi}_{\Delta t}(t)\right) \subset S^\gamma \left( \chi(t) \right) \quad \quad \forall t \in [0,T],\ \Delta t < \delta (\gamma). \label{922}
\end{align}

Now we fix an arbitrary function $b \in Y(\chi, T)$, hence there exists some $\gamma > 0$ such that $b$ is curl-free in $S^\gamma (\chi(t))$ for any $t \in [0,T]$. Then, by the second inclusion in (\ref{922}), $b$ is also curl-free in $S \left( \overline{\chi}_{\Delta t}(t)\right)$ for any $\Delta t < \delta (\gamma)$. In other words,
\begin{align}
b(t) \in W^k \quad \forall t \in ((k-1)\Delta t,k\Delta t],\ k=1,...,\frac{T}{\Delta t},\ \Delta t < \delta (\gamma), \label{794}
\end{align}

so we may use $b(t)$ as a test function for the discrete induction equation.\\
Next we take an arbitrary interval $I \subset (0,T)$ and an arbitrary open ball $U \subset \mathbb{R}^3$ such that $\overline{I\times U}\subset Q^S(\chi, T) \bigcap Q$. The first inclusion in (\ref{922}) implies
\begin{equation}
\text{curl}B = \lim_{\Delta t \rightarrow 0}\text{curl}\overline{B}_{\Delta t} = 0 \ \ \ \text{a.e. in } \overline{I \times U} \text{ and thus in } Q^S(\chi, T) \bigcap Q. \label{462}
\end{equation}

Now we take $b \in Y(\chi,T)$ and $\Delta t > 0$ sufficiently small such that (\ref{794}) holds true. For fixed $k$ we test the discrete induction equation (\ref{116}) by $b(t)$, $t \in ((k-1)\Delta t,k\Delta t]$, integrate over this interval and then sum up over all $k$ to see
\begin{align}
& \int _0^T\int _\Omega \partial_t B_{\Delta t} \cdot b\ dxdt = \sum_{k=1}^\frac{T}{\Delta t} \int_{(k-1)\Delta t}^{k\Delta t} \int _\Omega \frac{B_{\Delta t}^k - B_{\Delta t}^{k-1}}{\Delta t} \cdot b\ dxdt \nonumber \\
=& \int_0^T \int_\Omega \left(- \frac{1}{\sigma \mu} \text{curl} \overline{B}_{\Delta t} + \overline{u}_{\Delta t} \times \overline{B}'_{\Delta t} + \frac{1}{\sigma} \overline{J}_{\Delta t} - \frac{\epsilon}{\mu ^2}|\text{curl} \overline{B}_{\Delta t}|^2\text{curl} \overline{B}_{\Delta t} \right) \cdot \text{curl} b\ dxdt \nonumber \\
&-\int_0^T \int_\Omega \epsilon \text{curl} \left( \text{curl} \overline{B}_{\Delta t}\right) : \text{curl} \left( \text{curl} b \right)\ dxdt. \label{265}
\end{align}

An interpolation between $L^\infty(L^2)$ and $L^2(L^\infty)$ together with the estimates (\ref{29})--(\ref{344}) leads to the bounds
\begin{align}
\epsilon^\frac{1}{4}\left\| \overline{u}_{\Delta t} \right\|_{L^4(0,T;L^4(\Omega))} \leq c,\quad \epsilon^\frac{1}{4}\left\| \overline{B}_{\Delta t}' \right\|_{L^4(0,T;L^4(\Omega))} \leq c. \label{881}
\end{align}
Together with the bound (\ref{920}) this implies the existence of functions $z \in L^\frac{4}{3}((0,T)\times \Omega)$ and $z_1,z_2 \in L^2((0,T)\times \Omega)$ such that for chosen subsequences
\begin{align}
\epsilon \left| \text{curl} \overline{B}_{\Delta t} \right|^2 \text{curl} \overline{B}_{\Delta t} \rightharpoonup \epsilon z \quad &\text{in } L^\frac{4}{3}\left((0,T)\times \Omega\right), \label{851} \\
\overline{u}_{\Delta t} \times \overline{B}'_{\Delta t} \rightharpoonup z_1 \quad &\text{in } L^2\left((0,T)\times \Omega\right), \label{734} \\
\text{curl}\overline{B}'_{\Delta t} \times \overline{B}'_{\Delta t} \rightharpoonup z_2 \quad &\text{in } L^2\left((0,T)\times \Omega\right). \label{740}
\end{align}

With these convergences at hand we can pass to the limit in (\ref{265}) and obtain
\begin{align}
& - \int _0^T \int_\Omega B \cdot \partial_t b\ dxdt - \int_\Omega B_0 \cdot b(0,x)\ dx \nonumber \\
=& \int _0^T \int _\Omega \left[ -\frac{1}{\sigma \mu} \text{curl} B + z_1 + \frac{1}{\sigma} J - \frac{\epsilon}{\mu ^2} z \right] \cdot \text{curl} b - \epsilon \text{curl}\left( \text{curl} B\right): \text{curl}\left( \text{curl} b\right)\ dxdt \label{127}
\end{align}

for all $b\in Y(\chi,T)$. The main difficulty of this section is now to identify the limit functions $z_1$ and $z_2$. The limit function $z$ does not need to be identified, as it will vanish from the equation when we pass to the limit with $\epsilon \rightarrow 0$. We first note that
\begin{align}
z_1 \cdot \text{curl} b = 0 = (u \times B) \cdot \text{curl}b \quad &\text{a.e. in } Q^S(\chi,T)\bigcap Q, \label{735} \\
z_2 = 0 = \text{curl}B = \text{curl} B \times B \quad &\text{a.e. in } Q^S(\chi,T)\bigcap Q, \label{741}
\end{align}

where (\ref{735}) follows directly from $b \in Y(\chi,T)$ and (\ref{741}) follows in the same way as (\ref{462}). Hence it suffices to identify $z_1$ and $z_2$ in the fluid region. In order to do so, we choose an arbitrary interval $I = (a,d) \subset (0,T)$ and an arbitrary open ball $U \subset \Omega$ with $\overline{I \times U} \subset Q^F(\chi,T):= Q \setminus \overline{Q}^S(\chi,T)$. In (\ref{794}) we have seen that, for any sufficiently small $\Delta t >0$, functions from $\mathcal{D}(I \times U)$ are admissible test functions in (\ref{265}). By a density argument, (\ref{265}) may thus also be tested by any $b \in L^4(a,d;H_0^{2}(U))$, extended by $0$ outside of $(a, d)\times U$. This, together with the $L^4((0,T)\times \Omega)$-bound of $\overline{u}_{\Delta t}$ in (\ref{881}), leads to the dual estimate
\begin{equation}
\left\| \frac{\overline{B}_{\Delta t}(\cdot) - \overline{B}_{\Delta t}(\cdot - \Delta t) }{\Delta t} \right\|_{L^\frac{4}{3}(a, d;H^{-2}(U))} = \left\| \partial_t B_{\Delta t} \right\|_{L^\frac{4}{3}(a, d;H^{-2}(U))} \leq c. \label{266}
\end{equation}

From this estimate, we can now derive a corresponding estimate for the time-lagging interpolant $\overline{B}_{\Delta t}'$. Indeed, for arbitrary $b\in L^4 (a + \Delta t, d;H_0^{2,2}(U))$, we infer
\begin{align}
&\int_{a+\Delta t}^d\int_U \frac{\overline{B}'_{\Delta t}(t) - \overline{B}'_{\Delta t}(t-\Delta t)}{\Delta t} \cdot b(t)\ dxdt \nonumber \\
=&\int_{a}^{d-\Delta t}\int_U \frac{\overline{B}_{\Delta t}(t) - \overline{B}_{\Delta t}(t-\Delta t)}{\Delta t} \cdot b(t+\Delta t)\ dxdt \leq c \|b\|_{L^4\left(a+\Delta t,d;H_0^{2}(U)\right)}, \label{800}
\end{align}

so it holds
\begin{equation}
\left\| \frac{\overline{B}'_{\Delta t}(\cdot) - \overline{B}'_{\Delta t}(\cdot -\Delta t)}{\Delta t} \right\|_{L^\frac{4}{3}\left(a + \Delta t, d;H^{-2}(U)\right)} \leq c.\label{882}
\end{equation}

This gives us the conditions for the discrete Aubin-Lions Lemma \cite[Theorem 1]{dreherjungel}, which yields
\begin{equation}
\overline{B}'_{\Delta t} \rightarrow B \ \ \ \text{in } L^2(I;H^{-1}(U)). \label{125}
\end{equation}

By the arbitrary choice of $\overline{I \times U}$ this is sufficient to infer
\begin{align}
z_1 = u \times B \quad \text{a.e. in } Q^F(\chi,T),\quad z_2 = \text{curl}B \times B \quad \text{a.e. in } Q^F(\chi,T). \label{738}
\end{align}

\subsection{Continuity equation}

We test the discrete continuity equation (\ref{8}) by $\rho _{\Delta t}^k$, apply Young's inequality and sum over all $k=1,...,l$, $l \in \left\lbrace 1,...,\frac{T}{\Delta t} \right\rbrace$ which leads to
\begin{align}
\|\rho _{\Delta t}^l\|_{L^2(\Omega)}^2 + 2\Delta t\sum _{k=1}^{l}\epsilon \| \nabla \rho _{\Delta}^k \|_{L^2(\Omega)}^2 \leq \| \rho _0 \|_{L^2(\Omega)}^2 \quad \forall l \in \left\lbrace 1,...,\frac{T}{\Delta t} \right\rbrace. \label{921}
\end{align}

Hence, $\overline{\rho}_{\Delta t}$ is bounded in $L^2(0,T;H^{1}(\Omega))$ and we can find $\rho \in L^2(0,T;H^{1}(\Omega))$ such that for a chosen subsequence
\begin{equation}
\overline{\rho}_{\Delta t} \rightharpoonup \rho \quad \text{in } L^2\left(0,T;H^{1}(\Omega)\right). \label{853}
\end{equation}

Further, from the continuity equation we derive the dual estimates
\begin{equation}
\left\| \frac{\overline{\rho}_{\Delta t}(\cdot) - \overline{\rho}_{\Delta t}(\cdot - \Delta t)}{\Delta t} \right\|_{L^2\left(0,T;\left(H^{1}(\Omega)\right)^*\right)}\leq c,\quad \left\| \frac{\overline{\rho}'_{\Delta t}(\cdot) - \overline{\rho}'_{\Delta t}(\cdot -\Delta t)}{\Delta t} \right\|_{L^2\left(\Delta t,T;\left(H^{1}(\Omega)\right)^*\right)} \leq c \label{82}
\end{equation}

by the same arguments as the bounds (\ref{266}) and (\ref{882}) for the discrete time derivatives of $\overline{B}_{\Delta t}$ and $\overline{B}'_{\Delta t}$. In particular, we can again apply the discrete Aubin-Lions Lemma \cite[Theorem 1]{dreherjungel} to infer
\begin{equation}
\overline{\rho}_{\Delta t},\ \overline{\rho}'_{\Delta t} \rightarrow \rho \quad \text{in } L^q\left(0,T;L^q\left(\mathbb{R}^3\right)\right) \quad \forall 1\leq q < \infty, \quad \quad \underline{\rho} \leq \rho \leq \overline{\rho}\quad \text{a.e. in } [0,T]\times \mathbb{R}^3, \label{573}
\end{equation}

where the limit function $\rho$ has been extended by $\underline{\rho}$ outside of $\Omega$. Now, we sum the discrete equation (\ref{8}) over all $k=1,...,\frac{T}{\Delta t}$ and pass to the limit by means of (\ref{853}) and (\ref{573}). This yields
\begin{equation}
-\int_0^T \int _{\Omega} \rho \partial_t \psi \ dxdt - \int_{\Omega} \rho_0 \psi(0,x)\ dx = \int_0^T \int_{\Omega} (\rho u) \cdot \nabla \psi + \epsilon \rho \Delta \psi \ dxdt \quad \forall \psi \in \mathcal{D}([0,T)\times \Omega). \label{87}
\end{equation} 

Our next goal is to show strong convergence of $\nabla \overline{\rho}_{\Delta t}$, which is required for the limit passage in the momentum equation. The first bound in (\ref{82}) further implies that for a subsequence
\begin{equation}
\partial_t \rho_{\Delta t} \buildrel\ast\over\rightharpoonup \partial_t \rho \ \ \ \text{in } L^2\left(0,T;\left(H^{1}(\Omega)\right)^*\right). \label{84}
\end{equation} 

Consequently, the limit of the discrete continuity equation can also be expressed in the form
\begin{align}
\int_0^\tau \int _\Omega \partial_t \rho \psi - (\rho u) \cdot \nabla \psi + \epsilon \nabla \rho \cdot \nabla \psi \ dxdt = 0 \ \ \ \forall \psi \in L^2(0,T;H^{1}(\Omega)),\ \tau \in [0,T]. \label{576}
\end{align}

We now test (\ref{576}) by $\rho$ and compare it to the corresponding relation (\ref{921}) on the $\Delta t$--level, which will yield convergence of $\| \nabla \overline{\rho}_{\Delta t}\|_{L^2((0,\tau);L^2(\Omega))}$ and thus the desired strong convergence of $\nabla \overline{\rho}_{\Delta t}$. Indeed, testing (\ref{576}) by $\rho$ we obtain
\begin{equation}
\left\| \rho (\tau) \right\|_{L^2(\Omega)}^2 + 2\epsilon \int_0^\tau \int_\Omega |\nabla \rho|^2\ dxdt = \left\| \rho (0) \right\|_{L^2(\Omega)}^2. \label{578}
\end{equation}

Further, the inequality (\ref{921}) can be rewritten in the form
\begin{equation}
\|\overline{\rho} _{\Delta t}(l\Delta t - \gamma)\|_{L^2(\Omega)}^2 + 2\epsilon \int_{0}^{l\Delta t - \gamma} \int_\Omega |\nabla \overline{\rho}_{\Delta t}|^2\ dxdt \leq \| \rho_0 \|_{L^2(\Omega)}^2 \quad \forall \gamma \in [0,\Delta t),\ l \in \left\lbrace 1,...,\frac{T}{\Delta t} \right\rbrace. \label{888}
\end{equation}

Any $\tau \in (0,T]$ can be expressed as $\tau = l\Delta t - \gamma$ for some $l \in \left\lbrace 1,...,\frac{T}{\Delta t} \right\rbrace$ and $\gamma \in [0,\Delta t)$. Thus, subtracting (\ref{578}) from (\ref{888}) and making use of the strong convergence (\ref{573}), we infer that for a further subsequence
\begin{equation}
\lim_{\Delta t \rightarrow 0} \int_{0}^{\tau} \int_\Omega |\nabla \overline{\rho}_{\Delta t}|^2\ dxdt \leq \int_{0}^{\tau} \int_\Omega |\nabla \rho|^2\ dxdt \quad \text{for a.a. } \tau \in (0,T]. \label{579}
\end{equation}

On the other hand, for each such $\tau$ the weak lower semicontinuity of norms guarantees us the existence of $z_3 = z_3(\tau) \in \mathbb{R}$ such that for another subsequence
\begin{equation}
\left\| \nabla \overline{\rho}_{\Delta t} \right\|^2_{L^2(0,\tau;L^2(\Omega))} \rightarrow z_3^2 \geq \left\| \nabla \rho \right\|^2_{L^2(0,\tau;L^2(\Omega))}. \label{889}
\end{equation}

Combining (\ref{579}) and (\ref{889}), we infer that for almost all $\tau \in [0,T]$ there exists a subsequence for which
\begin{align}
\left\| \nabla \overline{\rho}_{\Delta t} \right\|_{L^2(0,\tau;L^2(\Omega))} \rightarrow \left\| \nabla \rho \right\|_{L^2(0,\tau;L^2(\Omega))}. \nonumber
\end{align}

In combination with the weak convergence (\ref{853}) and a diagonal argument, this implies the desired relation
\begin{equation}
\nabla \overline{\rho}_{\Delta t} \rightarrow \nabla \rho \ \ \ \text{in } L^2\left(0,\tau;L^2(\Omega)\right)\quad \text{for a.a. } \tau \in [0,T]. \label{697}
\end{equation}

Next, we show that the limit density satisfies a regularized and integrated version of the renormalized continuity equation (\ref{925}), which will be significant in the limit passage with respect to $\epsilon \rightarrow 0$. To this end we take an arbitrary smooth and convex function $\beta$ on $[\underline{\rho},\overline{\rho}]$ and test (\ref{8}) by $\beta'(\rho_{\Delta t}^k)$ for any $k=1,...,\frac{T}{\Delta t}$. By the convexity of $\beta$ and the fact that $\operatorname{div}\overline{u}_{\Delta t}' = 0$ this yields
\begin{align}
\int_0^\tau \int_\Omega \partial_t \rho_{\Delta t} \beta'(\overline{\rho}_{\Delta t})\ dx =& \Delta t\left( \sum_{k=1}^{l-1} \int_\Omega \frac{\rho_{\Delta t}^k - \rho_{\Delta t}^{k-1}}{\Delta t} \beta'(\rho_{\Delta t}^k)\ dx\right) + \gamma \int_\Omega \frac{\rho_{\Delta t}^l - \rho_{\Delta t}^{l-1}}{\Delta t} \beta'(\rho_{\Delta t}^l)\ dx \nonumber \\
=&-\int_0^\tau \int_\Omega \overline{u}_{\Delta t}' \cdot \nabla \overline{\rho}_{\Delta t} \beta '(\overline{\rho}_{\Delta t})\ dxdt - \int_0^\tau \int_\Omega \epsilon \nabla \overline{\rho}_{\Delta t} \nabla \beta'(\overline{\rho}_{\Delta t})\ dxdt \nonumber \\
=&-\int_0^\tau \int_\Omega \epsilon \left| \nabla \overline{\rho}_{\Delta t} \right|^2 \beta''(\overline{\rho}_{\Delta t})\ dxdt \leq 0 \label{894}
\end{align}

for any $\tau \in (0,T]$, and $l \in  \left\lbrace 1,...,\frac{T}{\Delta t} \right\rbrace$, $\gamma \in [0,\Delta t)$ chosen such that $\tau = l\Delta t - \gamma$. Since the derivatives of $\beta$ are bounded, the strong $L^2(H^{1})$-convergence of $\overline{\rho}_{\Delta t}$ (c.f.\ (\ref{573}), (\ref{697})) implies
\begin{align}
\beta''(\overline{\rho}_{\Delta t}) \buildrel\ast\over\rightharpoonup \beta''(\rho) \quad &\text{in } L^\infty((0,T)\times \Omega), \nonumber \\
\beta'(\overline{\rho}_{\Delta t}) \rightarrow \beta'(\rho) \quad &\text{in } L^2\left(0,\tau;H^{1}(\Omega)\right) \quad \text{for a.a. } \tau \in [0,T]. \nonumber
\end{align}

Using this in combination with (\ref{84}), we can pass to the limit in (\ref{894}) and obtain the desired relation
\begin{align}
\int_\Omega \beta(\rho (\tau))\ dx - \int_\Omega \beta(\rho_0)\ dx = \int_0^\tau \int_\Omega \partial_t \beta(\rho)\ dxdt = -\int_0^\tau \int_\Omega \epsilon \beta''(\rho)|\nabla \rho|^2\ dxdt \leq 0 \quad \text{for a.a. } \tau \in [0,T]. \label{895}
\end{align}

\subsection{Momentum equation} \label{Rothe limit in mom eq}

We test the discrete momentum equation (\ref{57}) by $\phi(t)$ for $\phi \in L^4(0,T;V_0^2(\Omega))$ and sum the result over all $k$. Using the Hölder inequality and the Gagliardo–Nirenberg interpolation inequality we estimate
\begin{align}
&\int _0^T \int_\Omega \epsilon \left( \nabla \overline{u}_{\Delta t} \nabla \overline{\rho}_{\Delta t}\right) \cdot \phi \ dx dt \nonumber \\
\leq  &\epsilon \| \nabla \overline{u}_{\Delta t} \|_{L^4(0,T;L^2(\Omega))} \| \nabla \overline{\rho}_{\Delta t} \|_{L^2(0,T;L^2(\Omega))}\| \phi \|_{L^4(0,T;L^\infty(\Omega))} \nonumber \\
\leq  &c\epsilon\left\| \overline{u}_{\Delta t} \right\|^\frac{1}{2}_{L^2(0,T;H^{2}(\mathbb{R}^3))}\left\| \overline{u}_{\Delta t} \right\|^\frac{1}{2}_{L^\infty(0,T;L^2(\mathbb{R}^3))} \| \nabla \overline{\rho}_{\Delta t} \|_{L^2(0,T;L^2(\Omega))}\| \phi \|_{L^4(0,T;L^\infty(\Omega))} \leq c. \nonumber
\end{align}

This allows us to infer the dual estimate
\begin{equation}
\left\| \frac{\overline{(\rho u)}_{\Delta t}(\cdot) - \overline{(\rho u)}_{\Delta t}(\cdot - \Delta t)}{\Delta t} \right\|_{L^\frac{4}{3}\left(0,T;\left(V_0^2(\Omega)\right)^*\right)} = \left\| \partial_t(\rho u)_{\Delta t} \right\|_{L^\frac{4}{3}\left(0,T;\left(V_0^2(\Omega)\right)^*\right)} \leq c. \label{73}
\end{equation}

We can now mimick the compactness results for the time-dependent incompressible Navier-Stokes equations, c.f.\ \cite[Theorem 2.4]{lions}. The estimate (\ref{73}) gives us the conditions for the discrete Aubin-Lions Lemma \cite[Theorem 1]{dreherjungel}, which allows us to deduce
\begin{equation}
 P\left(\overline{(\rho u)}_{\Delta t}\right) \rightarrow P(\rho u) \ \ \ \text{in } L^2\left(0,T;\left(V_0^2(\Omega)\right)^*\right), \nonumber
\end{equation}

where $P$ denotes the orthogonal projection of $L^2(\Omega)$ onto the space $V^0(\Omega)$ of weakly divergence-free $L^2$--functions. This, in combination with the $L^4((0,T)\times \Omega)$--bound (\ref{881}) of $\overline{u}_{\Delta t}$, leads to
\begin{equation}
\overline{u}_{\Delta t} \rightarrow u \ \ \ \text{in } L^q((0,T)\times \Omega) \quad \forall 1 \leq q < 4. \label{260}
\end{equation}

With this strong convergence at hand, we derive the following limit version of the momentum equation
\begin{align}
&- \int_0^T \int_\Omega \rho u \cdot \partial_t \phi \ dxdt - \int_\Omega \rho _0 u_0 \cdot \phi(0,x)\ dx \nonumber \\
=& \int _0^T \int _\Omega \rho (u \otimes u): \nabla \phi - 2\nu \mathbb{D}(u):\nabla \phi - \frac{1}{\eta} \rho \chi \left( u - \Pi \right) \cdot \phi \nonumber \\
&+ \rho g\cdot \phi + \frac{1}{\mu}\left( \text{curl}B \times B \right) \cdot \phi - \epsilon \left( \nabla u \nabla \rho \right) \cdot \phi - \epsilon \Delta u \cdot \Delta \phi \ dxdt \label{590}
\end{align}

for any $\phi \in \mathcal{D} ([0,T)\times \Omega)$ with $\text{div}\phi = 0$, where $\Pi$ was defined in (\ref{153}). Here we further used the strong convergence (\ref{697}) of $\nabla \overline{\rho}_{\Delta t}$ and the relations (\ref{740}), (\ref{741}), (\ref{738}) which identify the magnetic term in the limit equation. Now it only remains to identify $\Pi$. We start by remarking that
\begin{align}
\int_{\mathbb{R}^3} \rho (t) \chi(t) dx \geq \underline{\rho} |S| > 0 \quad \text{for a.a. } t \in [0,T]. \label{890}
\end{align}

We pick an arbitrary ball $B_R \subset \mathbb{R}^3$ with radius $R>0$, centered at $0$. The weak-$*$ convergence (\ref{796}) of $\overline{u}'_{\Delta t}$, the uniform convergence (\ref{919}) of the characteristic function and the strong convergence (\ref{573}) of the density yield that
\begin{align}
\int _{\mathbb{R}^3} \overline{\rho}_{\Delta t} \overline{\chi}'_{\Delta t} \ \overline{u}'_{\Delta t}\ dx &\buildrel\ast\over\rightharpoonup \int _{\mathbb{R}^3} \rho \chi \ u\ dx \quad \text{in } L^\infty\left((0,T)\times B_R\right), \nonumber \\
\int _{\mathbb{R}^3} \overline{\rho}_{\Delta t} \overline{\chi}'_{\Delta t}\ x\ dx &\rightarrow \int _{\mathbb{R}^3} \rho \chi \ x\ dx \quad \text{in } L^p\left((0,T)\times B_R\right) \quad \forall 1 \leq p < \infty, \nonumber \\
\int _{\mathbb{R}^3} \overline{\rho}_{\Delta t} \overline{\chi}'_{\Delta t}\ dx &\rightarrow \int _{\mathbb{R}^3} \rho \chi \ dx \quad \text{in } L^p\left((0,T)\times B_R\right) \quad \forall 1 \leq p < \infty. \nonumber
\end{align}

Combining the latter convergence with the bounds (\ref{216}), (\ref{890}) away from $0$ we further see
\begin{align}
\frac{1}{\int _{\mathbb{R}^3} \overline{\rho}_{\Delta t} \overline{\chi}'_{\Delta t}\ dx} &\rightarrow \frac{1}{\int _{\mathbb{R}^3} \rho \chi \ dx} \quad \text{in } L^p\left((0,T)\times B_R\right) \quad \forall 1\leq p < \infty, \nonumber \\
\frac{1}{\int _{\mathbb{R}^3} \overline{\rho}_{\Delta t} \overline{\chi}'_{\Delta t} \ dx} &\buildrel\ast\over\rightharpoonup \frac{1}{\int _{\mathbb{R}^3} \rho \chi \ dx} \quad \text{in } L^\infty\left((0,T)\times B_R\right) \nonumber
\end{align}

and altogether
\begin{equation}
\frac{\int_{\mathbb{R}^3} \overline{\rho}_{\Delta t} \overline{\chi}'_{\Delta t} \overline{u}'_{\Delta t}\ dx}{\int _{\mathbb{R}^3} \overline{\rho}_{\Delta t} \overline{\chi}'_{\Delta t}\ dx} = \overline{\left( u_G \right)}'_{\Delta t} \buildrel\ast\over\rightharpoonup \left(u_G\right)_{[\chi,\rho ,u]} := \frac{\int_{\mathbb{R}^3} \rho \chi u\ dx}{\int _{\mathbb{R}^3} \rho \chi \ dx}\ \ \ \text{in } L^\infty \left((0,T)\times B_R\right) \label{217}
\end{equation}

as well as
\begin{equation}
\frac{\int_{\mathbb{R}^3} \overline{\rho}_{\Delta t} \overline{\chi}'_{\Delta t} x\ dx}{\int _{\mathbb{R}^3} \overline{\rho}_{\Delta t} \overline{\chi}'_{\Delta t}\ dx} = \overline{a}'_{\Delta t} \rightarrow a_{[\chi,\rho]} := \frac{\int _{\mathbb{R}^3} \rho \chi x\ dx }{\int _{\mathbb{R}^3} \rho \chi \ dx}\ \ \ \text{in } L^p\left((0,T)\times B_R\right)\ \ \ \forall 1 \leq p < \infty. \label{203}
\end{equation}

Next, we consider the matrix
\begin{equation}
\overline{I}'_{\Delta t} = \int _{\mathbb{R}^3} \overline{\rho} _{\Delta t} \overline{\chi}'_{\Delta t} \left( | x - \overline{a}'_{\Delta t} |^2\text{id} - \left( x - \overline{a}'_{\Delta t} \right) \otimes \left( x - \overline{a}'_{\Delta t} \right) \right)\ dx, \nonumber
\end{equation}

for which similar arguments together with the already proven convergence (\ref{203}) lead to
\begin{equation}
\overline{I}'_{\Delta t} \rightarrow I_{[\chi,\rho]} := \int _{\mathbb{R}^3} \rho \chi \left( | x - a_{[\chi,\rho]} |^2\text{id} - \left( x - a_{[\chi,\rho]} \right) \otimes \left( x - a_{[\chi,\rho]} \right) \right)\ dx \ \ \ \text{in } L^p\left((0,T)\times B_R\right) \label{213}
\end{equation}

for any $1 \leq p < \infty$. From this and the bound of the eigenvalues of $\overline{I}_{\Delta t}'(t)$ away from $0$, c.f.\ (\ref{849}), it is possible to derive
\begin{equation}
(\overline{I}'_{\Delta t})^{-1} \rightarrow \left(I_{[\chi,\rho]}\right)^{-1}\ \ \ \text{in } L^p((0,T)\times B_R)\ \ \ \forall 1 \leq p < \infty. \label{214}
\end{equation}

This, together with arguments similar to the ones used for (\ref{217}), yields
\begin{align}
\overline{\omega}'_{\Delta t} =& \left(\overline{I}_{\Delta t}'\right)^{-1} \int_{\mathbb{R}^3} \overline{\rho}_{\Delta t} \overline{\chi}_{\Delta t}' \left( x - \overline{a}_{\Delta t}' \right) \times \overline{u}_{\Delta t}'\ dx \nonumber \\
\rightharpoonup \ \omega_{[\chi,\rho,u]} :=& I_{[\chi,\rho]}^{-1}\int_{\mathbb{R}^3} \rho \chi\left( \left( x - a_{[\chi,\rho]} \right) \times u \right)\ dx  \quad \text{in } L^p((0,T)\times B_R) \quad \forall 1 \leq p < \infty. \label{220}
\end{align}

Now (\ref{217}), (\ref{203}) and (\ref{220}) imply
\begin{equation}
\Pi = \left(u_G\right)_{[\chi,\rho,u]} + \omega _{[\chi,\rho, u]} \times (x-a_{[\chi,\rho]}) =: \Pi _{[\chi,\rho,u]}. \label{221}
\end{equation}

\subsection{Energy inequality} \label{Rothe limit in energy ineq}

In order to derive an energy inequality for the limit system, we first derive a slightly modified version of the discrete energy inequality (\ref{18}). More precisely, we again add the estimates (\ref{17}) and (\ref{340}) and sum over all $l=1,...,k$, $k \in \left\lbrace 1,...,\frac{T}{\Delta t} \right\rbrace$. Since each $\tau \in (0,T]$ can be written in the form $\tau = k\Delta t - \gamma$ for some $k\in \left\lbrace 1,...,\frac{T}{\Delta t} \right\rbrace$ and $\gamma \in [0,\Delta t)$, this leads to
\begin{align}
&\frac{1}{2}\left\| \sqrt{\overline{\rho}_{\Delta t}}(\tau) \overline{u}_{\Delta t}(\tau) \right\|_{L^2(\Omega)}^2 + \int_0^\tau \int_\Omega 2\nu \left| \nabla \overline{u}_{\Delta t}(t,x) \right|^2\ dxdt + \int_0^\tau \int_\Omega \epsilon \left| \Delta \overline{u}_{\Delta t}(t,x) \right|^2\ dxdt \nonumber \\
&+ \frac{1}{2\mu}\left\| \overline{B}_{\Delta t}(\tau) \right\|_{L^2(\Omega)}^2 + \int_0^\tau \int_\Omega \frac{\epsilon}{\mu ^3}\left| \text{curl} \overline{B}_{\Delta t}(t,x) \right|^4\ dxdt + \int_0^\tau \int_\Omega \frac{\epsilon}{\mu} \left| \Delta \overline{B}_{\Delta t}(t,x) \right|^2\ dxdt \nonumber \\
&+ \int_0^\tau \int_\Omega \frac{1}{\sigma \mu ^2} \left| \text{curl} \overline{B}_{\Delta t}(t,x) \right|^2\ dxdt \nonumber \\
\leq&\frac{1}{2}\left\| \sqrt{\rho}_0 u_0 \right\|_{L^2(\Omega)}^2 + \frac{1}{2\mu}\left\| B_0 \right\|_{L^2(\Omega)}^2 + \int_0^\tau \int_\Omega - \frac{1}{\eta} \overline{\rho}'_{\Delta t}(t,x) \overline{\chi}_{\Delta t}(t,x) \left( \overline{u}'_{\Delta t}(t,x) - \overline{\Pi}'_{\Delta t}(t,x) \right)\cdot \overline{u}_{\Delta t}(t,x) \nonumber \\
& + \overline{\rho}'_{\Delta t}(t,x) \overline{g}_{\Delta t}(t,x)\cdot \overline{u}_{\Delta t}(t,x) + \frac{1}{\mu} \left( \text{curl} \overline{B}'_{\Delta t}(t,x) \times \overline{B}'_{\Delta t}(t,x) \right) \cdot \overline{u}_{\Delta t}(t,x) \nonumber \\
& + \frac{1}{\mu} \left(\overline{u}_{\Delta t}(t,x) \times \overline{B}'_{\Delta t}(t,x) \right) \cdot \text{curl} \overline{B}_{\Delta t}(t,x) + \frac{1}{\sigma}\overline{J}_{\Delta t}(t,x) \cdot \text{curl} \overline{B}_{\Delta t}(t,x)\ dxdt + c\left[ \Delta t + (\Delta t)^\frac{1}{2} \right]. \nonumber
\end{align}

On the right-hand side of this inequality we can pass to the limit by using in particular the strong convergence (\ref{260}) of $\overline{u}_{\Delta t}$ and the relations (\ref{740}), (\ref{741}), (\ref{738}) which identify the limits of the mixed terms. Using the weak lower semicontinuity of norms on the left-hand side, we end up with
\begin{align}
&\int_\Omega \frac{1}{2} \rho (\tau) |u(\tau)|^2 + \frac{1}{2\mu}|B(\tau)|^2\ dx + \int_0^\tau \int_\Omega 2\nu \left| \nabla u(t,x) \right|^2 + \epsilon |\Delta u(t,x)|^2 + \frac{\epsilon}{\mu ^3}\left| z(t,x) \right|^\frac{4}{3} \nonumber \\
& + \frac{\epsilon}{\mu} \left| \Delta B(t,x) \right|^2 + \frac{1}{\sigma \mu^2} \left| \text{curl} B(t,x) \right|^2 + \frac{1}{\eta} \rho (t,x) \chi (t,x) \left| u(t,x) - \Pi _{[\chi,\rho ,u]}(t,x) \right|^2\ dxdt \nonumber \\
\leq& \int_\Omega \frac{1}{2} \rho_0 |u_0|^2 + \frac{1}{2}|B_0|^2\ dx + \int_0^\tau \int_\Omega \rho(t,x) g(t,x)\cdot u(t,x) + \frac{1}{\sigma}J(t,x) \cdot \text{curl} B(t,x)\ dxdt \label{339}
\end{align}

for almost all $\tau \in [0,T]$. Here, the mixed terms canceled each other by the identity
\begin{align}
\left( \operatorname{curl} B \times B \right) \cdot u = - \left( u \times B \right) \cdot \operatorname{curl} B, \nonumber
\end{align}

and the term involving $u - \Pi _{[\chi,\rho ,u]}$ was rewritten by means of the relation
\begin{equation}
\int_0^{\tau} \int_\Omega \rho(t,x) \chi(t,x) \left( u(t,x) - \Pi _{[\chi, \rho , u]}(t,x)\right) \cdot \Pi _{[\chi , \rho , u]}(t,x)\ dxdt = 0,\ \ \ \tau \in [0,T], \nonumber
\end{equation}

c.f. \cite[Lemma 3.1]{cottetmaitre}. In summary, we have shown

\begin{proposition}
\label{System on epsilon-level}

Let all the assumptions of Theorem \ref{mainresult} be satisfied and let $\epsilon > 0$. Assume in addition that
\begin{align}
\rho _0 \in H^{1}(\Omega),\quad \quad u_0,\ B_0 \in H^{2}(\Omega). \nonumber
\end{align}

Then, there exist
\begin{align}
\rho_\epsilon &\in \left\lbrace \psi \in L^2\left(0,T; H^{1}(\Omega) \right):\ \underline{\rho} \leq \psi \leq \overline{\rho}\ \text{a.e. in } Q \right\rbrace, \label{904} \\
\chi_\epsilon &\in C\left([0,T];L^p\left(\mathbb{R}^3\right)\right),\ 1 \leq p < \infty,\ \quad z_\epsilon \in L^\frac{4}{3}((0,T)\times \Omega), \label{945} \\
B_\epsilon &\in \bigg\lbrace b \in L^\infty \left(0,T;L^2(\Omega)\right) \bigcap L^2\left(0,T;H^{2}(\Omega)\right):\ \operatorname{div} b = 0\ \text{in } Q, \nonumber \\
&\quad \quad \operatorname{curl}b = 0 \ \text{in } Q^S(\chi_\epsilon,T) \bigcap Q,\ b \cdot \text{n}|_{\partial \Omega} = 0 \bigg\rbrace, \label{902} \\
u_\epsilon &\in  L^\infty \left(0,T;L^2(\Omega)\right) \bigcap L^2\left(0,T;V_0^2(\Omega)\right) \label{903}
\end{align}

such that
\begin{align}
-\int_0^T \int_{\mathbb{R}^3 } \chi_\epsilon \partial_t \Theta \ dxdt - \int _{\mathbb{R}^3 } \chi _0 \Theta (0,x)\ dx =& \int_0^T \int_{\mathbb{R}^3} \left( \chi_\epsilon \Pi_{[\chi_\epsilon, \rho_\epsilon, u_\epsilon]} \right) \cdot \nabla \Theta \ dxdt, \label{899} \\
-\int_0^T \int _{\Omega} \rho_\epsilon \partial_t \psi \ dxdt - \int_{\Omega} \rho_0 \psi(0,x)\ dx =& \int_0^T \int_{\Omega} (\rho_\epsilon u_\epsilon) \cdot \nabla \psi + \epsilon \rho_\epsilon \Delta \psi \ dxdt, \label{898} \\
- \int_0^T \int_\Omega \rho_\epsilon u_\epsilon \cdot \partial_t \phi \ dxdt - \int_\Omega \rho _0 u_0 \cdot \phi(0,x)\ dx =& \int _0^T \int _\Omega \rho_\epsilon (u_\epsilon \otimes u_\epsilon): \nabla \phi - 2\nu \mathbb{D}(u_\epsilon):\nabla \phi \nonumber \\
&- \frac{1}{\eta} \rho_\epsilon \chi_\epsilon \left( u_\epsilon - \Pi _{[\chi_\epsilon,\rho_\epsilon,u_\epsilon]} \right) \cdot \phi + \rho_\epsilon g\cdot \phi \nonumber \\
&+ \frac{1}{\mu}\left( \operatorname{curl}B_\epsilon \times B_\epsilon \right) \cdot \phi - \epsilon \left( \nabla u_\epsilon \nabla \rho_\epsilon \right) \cdot \phi \nonumber \\
&- \epsilon \Delta u_\epsilon \cdot \Delta \phi \ dxdt, \label{897} \\
- \int _0^T \int_\Omega B_\epsilon \cdot \partial_t b\ dxdt - \int_\Omega B_0 \cdot b(0,x)\ dx =& \int _0^T \int _\Omega \left[ -\frac{1}{\sigma \mu} \operatorname{curl} B_\epsilon + u_\epsilon \times B_\epsilon + \frac{1}{\sigma} J - \frac{\epsilon}{\mu ^2} z_\epsilon \right] \cdot \operatorname{curl} b \nonumber \\
&- \epsilon \operatorname{curl}\left( \operatorname{curl} B_\epsilon\right): \operatorname{curl}\left( \operatorname{curl} b\right)\ dxdt \label{900}
\end{align}

for all $\Theta \in \mathcal{D} ([0,T)\times \mathbb{R}^3)$, $\psi, \phi \in \mathcal{D} ([0,T)\times \Omega)$ and all $b\in Y(\chi,T)$. Moreover, these functions satisfy the energy inequality
\begin{align}
&\int_\Omega \frac{1}{2} \rho_\epsilon (\tau) |u_\epsilon(\tau)|^2 + \frac{1}{2\mu}|B_\epsilon(\tau)|^2\ dx + \int_0^\tau \int_\Omega 2\nu \left| \nabla u_\epsilon(t,x) \right|^2 + \epsilon |\Delta u_\epsilon(t,x)|^2 + \frac{\epsilon}{\mu ^3}\left| z_\epsilon(t,x) \right|^\frac{4}{3} \nonumber \\
& + \frac{\epsilon}{\mu} \left| \Delta B_\epsilon(t,x) \right|^2 + \frac{1}{\sigma \mu^2} \left| \operatorname{curl} B_\epsilon(t,x) \right|^2 + \frac{1}{\eta} \rho_\epsilon (t,x) \chi_\epsilon (t,x) \left| u_\epsilon(t,x) - \Pi _{[\chi_\epsilon,\rho_\epsilon ,u_\epsilon]}(t,x) \right|^2\ dxdt \nonumber \\
\leq& \int_\Omega \frac{1}{2} \rho_0 |u_0|^2 + \frac{1}{2}|B_0|^2\ dx + \int_0^\tau \int_\Omega \rho_\epsilon(t,x) g(t,x)\cdot u_\epsilon(t,x) + \frac{1}{\sigma}J(t,x) \cdot \operatorname{curl} B_\epsilon(t,x)\ dxdt \label{901}
\end{align}

for almost all $\tau \in [0,T]$ and the characteristic function $\chi_\epsilon$ is connected to the solution $X^{\Pi_{[\chi_\epsilon,\rho_\epsilon,u_\epsilon]}}$ of the initial value problem
\begin{align}
\frac{dX^{\Pi_{[\chi_\epsilon,\rho_\epsilon,u_\epsilon]}}(s;t,x)}{dt} &= \Pi_{[\chi_\epsilon,\rho_\epsilon,u_\epsilon]} \left(t, X^{\Pi_{[\chi_\epsilon,\rho_\epsilon,u_\epsilon]}} (s;t,x) \right), \quad \quad X^{\Pi_{[\chi_\epsilon,\rho_\epsilon,u_\epsilon]}}(s;s,x)= x \label{905}
\end{align}

by
\begin{align}\chi (t,x) = \chi_0\left(X^{\Pi_{[\chi_\epsilon,\rho_\epsilon,u_\epsilon]}}(t;0,x)\right). \label{906}
\end{align}

\end{proposition}

{\centering \section{Limit passage with respect to \texorpdfstring{$\epsilon \rightarrow 0$}{}} \label{regularizationlimit} \par
}

From the energy inequality (\ref{901}) we infer the existence of a constant $c>0$, independent of $\epsilon$, such that
\begin{align}
\|u_\epsilon\|_{L^\infty (0,T;L^2(\Omega))}+ \|B_\epsilon\|_{L^\infty (0,T;L^2(\Omega))} + \|u_\epsilon\|_{L^2\left(0,T;H^{1}(\Omega)\right)} + \| B_\epsilon\|_{L^2(0,T;H^{1}(\Omega))}&\leq c,  \label{353} \\
\epsilon^\frac{1}{2} \left\| \Delta u_\epsilon \right\|_{L^2((0,T)\times \Omega)} + \epsilon^\frac{3}{4} \left\| z_\epsilon \right\|_{L^{\frac{4}{3}}((0,T)\times \Omega)} + \epsilon^\frac{1}{2} \left\| \Delta B_\epsilon \right\|_{L^2((0,T)\times \Omega)} &\leq c. \label{354}
\end{align}

The continuity equation on the $\epsilon$-level tested by $\rho_\epsilon$, c.f.\ (\ref{578}), yields
\begin{align}
\epsilon \left\| \nabla \rho_\epsilon \right\|^2_{L^2((0,T)\times \Omega)} \leq c. \label{924}
\end{align}

Further, from the lower bound (\ref{890}) for the total mass of the solid we deduce, similarly to  (\ref{935}) and (\ref{849}), the estimates
\begin{align}
\left| a_{[\chi _{\epsilon},\rho _{\epsilon}]}(t) \right| \leq c,\quad \left| \left( u_G \right)_{[\chi _{\epsilon},\rho _{\epsilon},u_{\epsilon}]}(t) \right| \leq c \left\| u_\epsilon(t) \right\|_{L^2(\Omega)},\quad &\left| \omega_{[\chi_\epsilon, \rho_\epsilon, u_\epsilon]}(t) \right| \leq c \left\| u_\epsilon(t) \right\|_{L^2(\Omega)}, \label{927} \\
v \cdot \left( I_{[\chi_\epsilon, \rho_\epsilon]}(t) v\right) \geq c|v|^2 \ \ \forall v \in \mathbb{R}^3,& \label{439}
\end{align}

for the quantities $a_{[\chi _{\epsilon},\rho _{\epsilon}]}$, $\left( u_G \right)_{[\chi _{\epsilon},\rho _{\epsilon},u_{\epsilon}]}$, $\omega_{[\chi_\epsilon, \rho_\epsilon, u_\epsilon]}$ and $I_{[\chi_\epsilon, \rho_\epsilon]}$, introduced in (\ref{217})--(\ref{220}), with $c$ independent of $t$ and $\epsilon$ and therefore
\begin{equation}
\|\Pi _{[\chi _{\epsilon},\rho _{\epsilon},u_{\epsilon}]}(t)\|_{L^\infty(\Omega)} \leq c\|u_{\epsilon}(t)\|_{L^2(\Omega)} \ \ \ \text{for a.a. } t \in [0,T]. \label{437}
\end{equation}

By this, the bounds in (\ref{904}) for the density and the uniform bounds (\ref{353})--(\ref{924}) we find functions
\begin{align}
\rho &\in L^\infty ((0,T)\times  \Omega),\quad \Pi \in L^\infty \left(0,T;W^{1,\infty}_\text{loc} \left(\mathbb{R}^3\right)\right), \label{943} \\
B &\in \bigg\lbrace b \in L^\infty \left(0,T;L^2(\Omega)\right) \bigcap L^2\left(0,T;V^1(\Omega)\right):\ b \cdot \text{n}|_{\partial \Omega} = 0 \bigg\rbrace \label{946} \\
u &\in  L^\infty \left(0,T;L^2(\Omega)\right) \bigcap L^2\left(0,T;V_0^1(\Omega)\right) \label{947}
\end{align}

such that for chosen subsequences
\begin{align}
u_\epsilon \buildrel\ast\over\rightharpoonup u \ \ \ &\text{in } L^\infty \left(0,T;L^2(\Omega)\right), \ \ \ \ \ \ &u_\epsilon& \rightharpoonup u \ \ \ \ \text{in } L^2 \left(0,T;H^{1}(\Omega)\right), \label{355} \\
B_\epsilon \buildrel\ast\over\rightharpoonup B \ \ \ &\text{in } L^\infty \left(0,T;L^2(\Omega)\right), \ \ \ \ \ \ &B_\epsilon& \rightharpoonup B \ \ \ \text{in } L^2 \left(0,T;H^{1}(\Omega)\right), \label{356} \\
\rho_\epsilon \buildrel\ast\over\rightharpoonup \rho  \ \ \ &\text{in } L^\infty \left(0,T;L^\infty(\Omega)\right), \ \ \ \ \ \  &\Pi_{[\chi_\epsilon , \rho _\epsilon , u_\epsilon]}& \buildrel\ast\over\rightharpoonup \Pi \ \ \ \text{in } L^\infty \left(0,T;W^{1,\infty}_\text{loc} \left(\mathbb{R}^3\right)\right) \label{357}
\end{align}

and
\begin{align}
\epsilon \nabla \rho_\epsilon,\ \epsilon \Delta u_\epsilon,\ \epsilon \Delta B_\epsilon \rightarrow 0 \quad \text{in } L^2((0,T)\times \Omega),\quad \quad \epsilon z_\epsilon \rightarrow 0 \quad \text{in } L^{\frac{4}{3}}((0,T)\times \Omega). \label{923}
\end{align}

\subsection{Characteristic function} \label{characteristicfunctionregularizationlimit}

The transport equation (\ref{899}), the equation (\ref{905}) for the associated characteristics and the estimates (\ref{927}) correspond directly to the conditions for Lemma \ref{strongconvergencecharacteristics} and Remark \ref{strongconvergencecharacteristicsremark}, which therefore yield
\begin{align}
X_\epsilon^{\Pi_{[\chi_\epsilon , \rho _\epsilon , u_\epsilon]}} \rightarrow X^{\Pi} \quad &\text{in } C([0,T]\times [0,T];C_\text{loc}\left(\mathbb{R}^3 )\right), \label{670} \\
\chi_\epsilon \rightarrow \chi \quad &\text{in } C\left([0,T];L^p\left(\mathbb{R}^3 \right)\right) \quad \forall 1 \leq p < \infty,\quad \quad \chi (t, x) = \chi_0\left(X^{\Pi}(t;0,x)\right), \label{492}
\end{align}

where $\chi$ and $X^{\Pi}$ are the unique solutions to
\begin{align}
-\int_0^T \int_{\mathbb{R}^3 } \chi \partial_t \Theta dxdt - \int _{\mathbb{R}^3 } \chi _0 \Theta (0,x)\ dx &= \int_0^T \int_{\mathbb{R}^3} \left(\chi  \Pi \right) \cdot \nabla \Theta \ dxdt \quad \forall \Theta \in \mathcal{D} \left([0,T)\times \mathbb{R}^3\right), \label{493} \\
\frac{dX^{\Pi}(s;t,x)}{dt} &= \Pi \left(t, X^{\Pi} (s;t,x) \right), \quad \quad X^{\Pi}(s;s,x) = x \label{671}.
\end{align}

\subsection{Induction equation} \label{regularizationlimitinduction}

Interpolating the bounds for $B_\epsilon$ in $L^\infty(0,T;L^2(\Omega))$ and $L^2(0,T;L^6(\Omega))$ we see that $B_\epsilon$ is also bounded in $L^3((0,T)\times \Omega)$. Hence, using the Hölder inequality, we find $z_4,z_5 \in L^\frac{6}{5}((0,T)\times \Omega)$ such that for selected subsequences
\begin{align}
u_\epsilon \times B_\epsilon \rightharpoonup z_4 \quad \text{in } L^\frac{6}{5}\left((0,T)\times \Omega\right), \quad \quad \text{curl}B_\epsilon \times B_\epsilon \rightharpoonup z_5 \quad \text{in } L^\frac{6}{5}\left((0,T)\times \Omega\right). \label{804}
\end{align}

Further, for any $\gamma > 0$ we again find, by (\ref{670}) and (\ref{492}), some $\delta (\gamma)>0$ such that
\begin{equation}
S_\gamma (\chi (t)) \subset S(\chi_\epsilon (t)) \subset S^{\frac{\gamma}{2}} (\chi (t)) \subset S^\gamma (\chi (t)) \quad \forall t \in [0,T],\ \epsilon < \delta (\gamma).\label{860}
\end{equation}

We fix arbitrary $b \in Y(\chi,T)$, so $b$ is curl-free in $S^\gamma (\chi (t))$ for some $\gamma > 0$ and all $t \in [0,T]$. Now (\ref{860}) implies that $b$ is also curl-free in a $\frac{\gamma}{2}$-neighbourhood of the solid region on the $\epsilon$-level for all sufficiently small $\epsilon>0$. In particular it holds $b \in Y(\chi _\epsilon,T)$ for all such $\epsilon$. Thus, letting $\epsilon \rightarrow 0$ in (\ref{900}), we obtain
\begin{align}
& - \int _0^{T} \int_\Omega B \cdot \partial_t b\ dxdt - \int_\Omega B_0 \cdot b(0)\ dx = \int _0^{T} \int _\Omega \left[ -\frac{1}{\sigma \mu} \text{curl} B + z_4 + \frac{1}{\sigma} J \right] \cdot \text{curl} b\ dxdt \label{515}
\end{align}

for any $b \in Y(\chi,T)$, where the regularization terms vanished as stated in (\ref{923}). It remains to identify $z_4$ and $z_5$. On the solid domain, we can argue as in Section \ref{inductionrothelimit} and see from the fact that $B_\epsilon$ is curl-free in $Q^S(\chi_\epsilon,T) \bigcap Q$ and (\ref{860}) that
\begin{align}
z_4 \cdot \text{curl}b = 0 = (u \times B) \cdot \text{curl} b, \quad \text{curl}B \times B = \text{curl} B = 0 = z_5 \quad \quad \text{a.e. in } Q^S(\chi,T) \bigcap Q \label{862}
\end{align}

for $b \in Y(\chi,T)$. In the fluid region we again consider an arbitrary set of the form $\overline{I \times U}\subset Q^F(\chi,T)$, where $I\subset (0,T)$ is an interval and $U\subset \Omega$ is a ball. For any sufficiently small $\epsilon > 0$ the first inclusion in (\ref{860}) implies that, for all functions $\psi \in \mathcal{D}(I)$ and $b \in \mathcal{D}(U)$ extended by $0$ outside of $I$ and $U$, the product $\psi b$ is an admissible test function in the induction equation (\ref{900}) on the $\epsilon$-level. This, together with the uniform estimates (\ref{353}), (\ref{354}) leads to the dual estimate
\begin{align}
\left\| \partial_t \int_U B_\epsilon \cdot b dx \right\|_{L^\frac{4}{3}(I)} \leq c. \label{884}
\end{align}

This allows us to apply the Arzelà-Ascoli theorem and deduce
\begin{equation}
B_\epsilon \rightarrow B \quad \text{in } C_ \text{weak}\left(\overline{I};L^2(U)\right)\ \text{and thus in } L^p\left(I;H^{-1}(U)\right)\quad \forall 1 \leq p < \infty. \label{393}
\end{equation}

Hence, writing
\begin{align}
\int_I \int_U \left(\text{curl} B_\epsilon \times B_\epsilon\right) \cdot b\ dxdt =& \int_I \int_U \text{div}\left( B_\epsilon \otimes B_\epsilon \right) \cdot b - \nabla \left( \frac{1}{2} \left| B_\epsilon \right|^2 \right) \cdot b\ dxdt, \quad b \in \mathcal{D}(I \times U), \label{394}
\end{align}

which allows us, after integration by parts, to shift the derivatives to the test function $b$, we conclude the desired identities
\begin{equation}
z_4 = u \times B \quad \text{a.e. in } Q^F(\chi,T),\quad \quad z_5 = \text{curl} B \times B \quad \text{a.e. in } Q^F(\chi,T).  \label{397}
\end{equation}

\subsection{Continuity equation} \label{conteqregularizationlimit}

We test the continuity equation (\ref{898}) by $\psi \Phi$, where $\psi \in \mathcal{D}(0,T)$ and $\Phi \in \mathcal{D}(\Omega)$, to find that
\begin{align}
\left\| \partial_t \int_\Omega \rho_\epsilon \Phi dx \right\|_{L^2(0,T)} \leq c. \nonumber
\end{align}

This again gives us the conditions for the Arzelà-Ascoli theorem, from which we obtain
\begin{align}
\rho _\epsilon \rightarrow \rho \quad \text{in } C_{\text{weak}}\left( [0,T];L^2(\Omega) \right)\ \text{and thus in } L^p\left(0,T;\left(H^{1}(\Omega)\right)^*\right)\quad \forall 1 \leq p < \infty. \nonumber
\end{align}

Combining this with the weak convergence (\ref{355}) of $u_\epsilon$ and the fact that $\epsilon \nabla \rho_\epsilon$ converges to $0$ in $L^2((0,T)\times \Omega)$ according to (\ref{923}), we may pass to the limit in (\ref{898}) and obtain
\begin{equation}
-\int_0^T \int _{\Omega} \rho \partial_t \psi dxdt - \int_{\Omega} \rho_0 \psi(0,x )\ dx = \int_0^T \int_\Omega (\rho u) \cdot \nabla \psi \ dxdt \quad \forall \psi \in \mathcal{D}([0,T)\times \Omega). \label{436}
\end{equation}

Since $\rho \in L^2((0,T)\times \Omega)$, the transport theorem by DiPerna and Lions \cite{dipernalions} implies that $\rho$ also satisfies, in the sense of distributions, the renormalized continuity equation (\ref{925}) for any bounded $\beta \in C^1(\mathbb{R})$ vanishing near $0$ and such that also $(\beta '(1 + |\cdot|))^{-1}$ is bounded. As $\rho$ is bounded from above and away from $0$, we can actually choose $\beta (z) = z \text{ln}(z)$. Using the same choice in the corresponding relation (\ref{895}) on the $\epsilon$-level, letting $\epsilon \rightarrow 0$ and comparing the results, we conclude
\begin{align}
\lim_{\epsilon \rightarrow 0} \int_\Omega \rho_\epsilon \text{ln}(\rho_\epsilon (\tau))\ dx \leq \int_\Omega \rho \text{ln}(\rho (\tau))\ dx \quad \text{for a.a. } \tau \in [0,T]. \nonumber
\end{align}

Following e.g.\ \cite[Theorem 10.20]{singularlimits}, this implies, by the strict convexity of $z \mapsto z \text{ln}(z)$, that
\begin{align}
\rho_\epsilon \rightarrow \rho \quad \text{a.e. in } (0,T)\times \Omega. \nonumber
\end{align}

In particular it follows
\begin{align}
\rho_\epsilon \rightarrow \rho \quad \text{in } L^p \left((0,T)\times \mathbb{R}^3\right)\quad \forall 1\leq p <\infty,\quad \quad \underline{\rho} \leq \rho \leq \overline{\rho}\quad \text{a.e. in } [0,T] \times \mathbb{R}^3, \label{487}
\end{align}

where $\rho$ has again been extended by $\underline{\rho}$ outside of $\Omega$.

\subsection{Momentum equation} \label{plaplacianlimitmomentum}

In order to pass to the limit in the momentum equation we further need strong convergence of the velocity field. We test the momentum equation (\ref{897}) on the $\epsilon$-level by $\psi \Phi$, where $\psi \in \mathcal{D}(0,T)$ and $\Phi \in \mathcal{D}(\Omega)$ with $\text{div}\Phi = 0$. This yields
\begin{align}
\left\| \partial _t \int_\Omega P(\rho _\epsilon u_\epsilon) \cdot \Phi dx \right\|_{L^\frac{4}{3}(0,T)} \leq c, \label{931}
\end{align}

where $P$ again denotes the orthogonal projection of $L^2(\Omega)$ onto $V^0(\Omega)$. The estimate (\ref{931}) leads, under exploitation of the Arzelà-Ascoli theorem, to
\begin{align}
P(\rho _\epsilon u_\epsilon) \rightarrow P(\rho u) \quad \text{in } C_\text{weak}([0,T];L^2(\Omega))\ \text{and thus in } L^2\left(0,T;\left( V_0^1(\Omega) \right)^*\right). \nonumber
\end{align}

By the same arguments as in the proof of the classical compactness results for the incompressible Navier-Stokes equations, c.f.\ \cite[Theorem 2.4]{lions}, this yields strong convergence of $u_\epsilon$ in $L^2((0,T)\times \Omega)$ and in particular
\begin{equation}
\rho _\epsilon u_\epsilon \otimes u_\epsilon \rightharpoonup \rho u \otimes u \quad \text{in } L^2\left(0,T;L^\frac{3}{2}(\Omega)\right). \label{572}
\end{equation}

Moreover, we can use the strong convergence (\ref{492}) of the characteristic function and the strong convergence (\ref{487}) of the density to identify the limit function $\Pi$ from (\ref{357}) as $\Pi = \Pi_{[\chi,\rho,u]}$ just as in (\ref{221}). Combining this with (\ref{804}), (\ref{862}), (\ref{397}) for the identification of the magnetic term and (\ref{572}), we can pass to the limit in (\ref{897}). The regularization terms again vanish as stated in (\ref{923}) and so we end up with
\begin{align}
&- \int _0^{T} \int_\Omega \rho u \cdot \partial_t \phi \ dxdt - \int_\Omega \rho _0 u_0 \cdot \phi(0,x)\ dx\nonumber \\
=& \int _0^{T} \int _\Omega \rho (u \otimes u): \nabla \phi - 2\nu \mathbb{D}(u):\nabla \phi - \frac{1}{\eta} \rho \chi \left( u - \Pi _{[\chi,\rho,u]} \right) \cdot \phi + \rho g\cdot \phi + \frac{1}{\mu}\left( \text{curl}B \times B \right) \cdot \phi \ dxdt \label{465}
\end{align}

for any $\phi \in \mathcal{D}([0,T)\times \Omega)$ with $\text{div}\phi = 0$.

\subsection{Energy inequality}

We drop the (nonnegative) regularization terms from the left-hand side of the energy inequality (\ref{901}). Using weak lower semicontinuity of norms, we then let $\epsilon$ tend to $0$ and obtain
\begin{align}
&\int_\Omega \frac{1}{2} \rho (\tau) |u(\tau)|^2 + \frac{1}{2}|B(\tau)|^2\ dx + \int_0^\tau \int_\Omega 2\nu \left| \nabla u(t,x) \right|^2 + \frac{1}{\sigma \mu^2} \left| \text{curl} B(t,x) \right|^2 \nonumber \\
&+ \frac{1}{\eta} \rho (t,x) \chi (t,x) \left| \left( u(t,x) - \Pi _{[\chi,\rho,u]}(t,x) \right) \right|^2\ dxdt \nonumber \\
\leq& \int_\Omega \frac{1}{2} \rho (0) |u(0)|^2 + \frac{1}{2}|B(0)|^2\ dx + \int_0^\tau \int_\Omega \rho (t,x) g(t,x)\cdot u(t,x) + \frac{1}{\sigma}J(t,x) \cdot \text{curl} B(t,x)\ dxdt  \label{478}
\end{align}

for almost all $\tau \in [0,T]$. Altogether we have shown

\begin{proposition}
\label{System on eta-level}

Let all the assumptions of Theorem \ref{mainresult} be satisfied and let $\eta > 0$. Assume in addition that
\begin{align}
\rho _{0,\eta} \in H^{1}(\Omega),\quad \quad u_{0,\eta},\ B_{0,\eta} \in H^{2}(\Omega). \nonumber
\end{align}

Then there exist
\begin{align}
\rho_\eta &\in \left\lbrace \psi \in L^\infty \left((0,T) \times \Omega \right):\ \underline{\rho} \leq \psi \leq \overline{\rho}\ \text{a.e. in } Q \right\rbrace, \label{907} \\
\chi_\eta &\in C\left([0,T];L^p\left(\mathbb{R}^3\right)\right),\ 1 \leq p < \infty, \label{908} \\
B_\eta &\in \left\lbrace b \in L^\infty \left(0,T;L^2(\Omega)\right) \bigcap L^2\left(0,T;H^{1}(\Omega)\right):\ \operatorname{div} b = 0\ \text{in } Q,\right. \nonumber \\
&\quad \quad \left. \operatorname{curl}b = 0 \ \text{in } Q^S(\chi_\eta,T) \bigcap Q,\ b \cdot \text{n}|_{\partial \Omega} = 0 \right\rbrace, \label{909} \\
u_\eta &\in  L^\infty \left(0,T;L^2(\Omega)\right) \bigcap L^2\left(0,T;V_0^1(\Omega)\right) \label{910}
\end{align}

such that
\begin{align}
-\int_0^T \int_{\mathbb{R}^3 } \chi_\eta \partial_t \Theta dxdt - \int _{\mathbb{R}^3 } \chi _0 \Theta (0,x)\ dx =& \int_0^T \int_{\mathbb{R}^3} \left(\chi_\eta \Pi_{[\chi_\eta, \rho_\eta, u_\eta]} \right) \cdot \nabla \Theta \ dxdt, \label{911} \\
-\int_0^T \int _{\Omega} \rho_\eta \partial_t \psi dxdt - \int_{\Omega} \rho_{0,\eta} \psi(0,x)\ dx =& \int_0^T \int_\Omega \left(\rho _\eta u_\eta \right) \cdot \nabla \psi \ dxdt, \label{912}\\
- \int_0^T \int_\Omega \rho_\eta u_\eta \cdot \partial_t \phi \ dxdt - \int_\Omega \rho _{0,\eta} u_{0,\eta} \cdot \phi(0,x)\ dx =& \int _0^T \int _\Omega \rho_\eta (u_\eta \otimes u_\eta): \nabla \phi - 2\nu \mathbb{D}(u_\eta):\nabla \phi \nonumber \\
&- \frac{1}{\eta} \rho_\eta \chi_\eta \left( u_\eta - \Pi _{[\chi_\eta,\rho_\eta,u_\eta]} \right) \cdot \phi + \rho_\eta g\cdot \phi \nonumber \\
&+ \frac{1}{\mu}\left( \operatorname{curl}B_\eta \times B_\eta \right) \cdot \phi \ dxdt, \label{913} \\
- \int _0^T \int_\Omega B_\eta \cdot \partial_t b\ dxdt - \int_\Omega B_{0,\eta} \cdot b(0,x)\ dx =& \int _0^T \int _\Omega \left[ -\frac{1}{\sigma \mu} \operatorname{curl} B_\eta + u_\eta \times B_\eta + \frac{1}{\sigma} J \right] \cdot \operatorname{curl} b\ dxdt \label{914}
\end{align}

for all $\Theta \in \mathcal{D} ([0,T)\times \mathbb{R}^3),$ $\psi, \phi \in \mathcal{D} ([0,T)\times \Omega)$ and all $b\in Y(\chi,T)$. Moreover, these functions satisfy the energy inequality
\begin{align}
&\int_\Omega \frac{1}{2} \rho_\eta (\tau) |u_\eta(\tau)|^2 + \frac{1}{2\mu}|B_\eta(\tau)|^2\ dx + \int_0^\tau \int_\Omega 2\nu \left| \nabla u_\eta(t,x) \right|^2 \nonumber \\
& + \frac{1}{\sigma \mu^2} \left| \operatorname{curl} B_\eta(t,x) \right|^2 + \frac{1}{\eta} \rho_\eta (t,x) \chi_\eta (t,x) \left| u_\eta(t,x) - \Pi _{[\chi_\eta,\rho_\eta ,u_\eta]}(t,x) \right|^2\ dxdt \nonumber \\
\leq& \int_\Omega \frac{1}{2} \rho_{0,\eta} |u_{0,\eta}|^2 + \frac{1}{2}|B_{0,\eta}|^2\ dx + \int_0^\tau \int_\Omega \rho_\eta(t,x) g(t,x)\cdot u_\eta(t,x) + \frac{1}{\sigma}J(t,x) \cdot \text{curl} B_\eta(t,x)\ dxdt \label{915}
\end{align}

for almost all $\tau \in [0,T]$ and the characteristic function $\chi_\eta$ is connected to the solution $X^{\Pi_{[\chi_\eta,\rho_\eta,u_\eta]}}$ of the initial value problem
\begin{align}
\frac{dX^{\Pi_{[\chi_\eta,\rho_\eta,u_\eta]}}(s;t,x)}{dt} &= \Pi_{[\chi_\eta,\rho_\eta,u_\eta]} \left(t, X^{\Pi_{[\chi_\eta,\rho_\eta,u_\eta]}} (s;t,x) \right), \quad \quad X^{\Pi_{[\chi_\eta,\rho_\eta,u_\eta]}}(s;s,x)= x \label{916}
\end{align}

by
\begin{align}\chi (t,x) = \chi_0\left(X^{\Pi_{[\chi_\eta,\rho_\eta,u_\eta]}}(t;0,x)\right). \label{917}
\end{align}

\end{proposition}

{\centering \section{Limit passage with respect to \texorpdfstring{$\eta \rightarrow 0$}{}} \label{projectionlimit} \par
}

\subsection{Uniform bounds and convergent terms} \label{etabounds}

In order to prove Theorem \ref{mainresult}, we assume in this section further that the regularized initial data we had chosen on the $\Delta t$-level and the $\epsilon$-level satisfy       
\begin{align}
\rho_{0,\eta} \rightarrow \rho_0 \quad \text{in } L^2(\Omega), \quad \quad u_{0,\eta} \rightarrow u_0 \quad \text{in } L^2(\Omega), \quad \quad B_{0,\eta} \rightarrow B_0 \quad \text{in } L^2(\Omega), \label{869}
\end{align}

where $\rho_0$, $u_0$, $B_0$ denote the initial data in Theorem \ref{mainresult}. The energy inequality (\ref{915}) implies the existence of a constant $c>0$, independent of $\eta$, such that
\begin{align}
\|u_\eta\|_{L^\infty (0,T;L^2(\Omega))}+ \|B_\eta\|_{L^\infty (0,T;L^2(\Omega))} + \|u_\eta\|_{L^2\left(0,T;H^{1}(\Omega)\right)} + \| B_\eta\|_{L^2(0,T;H^{1}(\Omega))}\leq& c, \label{930} \\
\frac{1}{\eta^\frac{1}{2}} \left\| \chi_\eta \left( u_\eta - \Pi _{[\chi_\eta,\rho_\eta ,u_\eta]} \right) \right\|_{L^2(0,T;L^2(\Omega))} \leq&c, \label{932}
\end{align}

and as in the corresponding estimates (\ref{927})--(\ref{437}) on the $\epsilon$-level, we deduce that
\begin{align}
\left| a_{[\chi _{\eta},\rho _{\eta}]}(t) \right| \leq c,\quad \left| \left( u_G \right)_{[\chi _{\eta},\rho _{\eta},u_{\eta}]}(t) \right| &\leq c \left\| u_\eta(t) \right\|_{L^2(\Omega)},\quad \left| \omega_{[\chi_\eta, \rho_\eta, u_\eta]}(t) \right| \leq c \left\| u_\eta(t) \right\|_{L^2(\Omega)}, \label{863} \\
v \cdot \left( I_{[\chi_\eta, \rho_\eta]}(t)v\right) &\geq c|v|^2 \quad \forall v \in \mathbb{R}^3, \label{928} \\
\left\|\Pi _{[\chi _{\eta},\rho _{\eta},u_{\eta}]}(t,\cdot)\right\|_{W^{1,\infty}(\Omega)} &\leq c \left\| u_{\eta} (t) \right\|_{L^2(\Omega)}\quad \text{for} \text{ a.a. } t \in [0,T] \label{648}
\end{align}

with $c$ independent of $\eta$ and $t$. The above bounds, together with the uniform bounds for the density in (\ref{907}), allow us to find functions
\begin{align}
\rho &\in L^\infty ((0,T)\times  \Omega), \label{942} \\
B &\in \bigg\lbrace b \in L^\infty \left(0,T;L^2(\Omega)\right) \bigcap L^2\left(0,T;V^1(\Omega)\right):\ b \cdot \text{n}|_{\partial \Omega} = 0 \bigg\rbrace, \label{891} \\
u &\in  L^\infty \left(0,T;L^2(\Omega)\right) \bigcap L^2\left(0,T;V_0^1(\Omega)\right) \label{892}
\end{align}

such that for extracted subsequences
\begin{align}
u_\eta \buildrel\ast\over\rightharpoonup u \ \ \ &\text{in } L^\infty \left(0,T;L^2(\Omega)\right), \ \ \ \ \ \ &u_\eta& \rightharpoonup u \ \ \ \ \ \ \ \ \ \ \ \text{in } L^2 \left(0,T;H^{1}(\Omega)\right), \label{581} \\
B_\eta \buildrel\ast\over\rightharpoonup B \ \ \ &\text{in } L^\infty \left(0,T;L^2(\Omega)\right), \ \ \ \ \ \ &B_\eta& \rightharpoonup B \ \ \ \ \ \ \ \ \ \ \text{in } L^2 \left(0,T;H^{1}(\Omega)\right), \label{582} \\
\rho_\eta \buildrel\ast\over\rightharpoonup \rho \ \ \ \ &\text{in } L^\infty \left(0,T;L^\infty(\Omega)\right), &\Pi_{[\chi_\eta , \rho _\eta , u_\eta]}& \buildrel\ast\over\rightharpoonup \Pi_{[\chi,\rho,u]} \ \ \ \text{in } L^\infty \left(0,T;W^{1,\infty}_\text{loc} \left(\mathbb{R}^3\right)\right). \label{583}
\end{align}

The identification of the limit function $\Pi_{[\chi,\rho,u]} = \left(u_G\right)_{[\chi,\rho,u]} + \omega _{[\chi,\rho, u]} \times (x-a_{[\chi,\rho]})$ in (\ref{583}) can be obtained as in the derivation of (\ref{221}), under the exploitation of the strong convergence of $\chi_\eta$ in (\ref{867}) and $\rho_\eta$ in (\ref{618}) below.

\subsection{Characteristic function}

The transport equation (\ref{911}), the equation (\ref{916}) for the corresponding characteristics and the bounds (\ref{863}) allow us to once more apply Lemma \ref{strongconvergencecharacteristics} and Remark \ref{strongconvergencecharacteristicsremark}, which yield
\begin{align}
X_\eta^{\Pi_{[\chi_\eta , \rho _\eta , u_\eta]}} \rightarrow X^{\Pi_{[\chi,\rho,u]}} \quad &\text{in } C([0,T]\times [0,T];C_\text{loc}\left(\mathbb{R}^3 )\right), \label{865} \\
\chi_\eta \rightarrow \chi \quad &\text{in } C\left([0,T];L^p\left(\mathbb{R}^3 \right)\right) \quad \forall 1\leq p < \infty,\quad \quad \chi (t, x) = \chi_0\left(X^{\Pi}(t;0,x)\right), \label{867}
\end{align}

where $X^{\Pi_{[\chi,\rho,u]}}$ and $\chi$ denote the unique solutions of
\begin{align}
-\int_0^T \int_{\mathbb{R}^3 } \chi \partial_t \Theta dxdt - \int _{\mathbb{R}^3 } \chi _0 \Theta (0,x)\ dx =& \int_0^T \int_{\mathbb{R}^3} \left( \chi \Pi_{[\chi,\rho,u]} \right) \cdot \nabla \Theta \ dxdt \quad \forall \Theta \in \mathcal{D} \left([0,T)\times \mathbb{R}^3\right), \label{868} \\
\frac{dX^{\Pi_{[\chi,\rho,u]}}(s;t,x)}{dt} =& \Pi_{[\chi,\rho,u]} \left(t, X^{\Pi_{[\chi,\rho,u]}} (s;t,x) \right), \quad \quad X^{\Pi_{[\chi,\rho,u]}}(s;s,x) = x. \label{866}
\end{align}

\subsection{Continuity equation}
 
For the strong convergence of the density we can apply classical compactness results for the incompressible Navier-Stokes equations, c.f.\ \cite[Theorem 2.4, Remark 2.4 3)]{lions}, and infer that
\begin{align}
\rho_\eta \rightarrow \rho \quad \text{in } C\left([0,T];L^p\left(\mathbb{R}^3\right)\right)\quad \forall 1 \leq p < \infty, \label{618}
\end{align}

with $\rho$ once again extended by $\underline{\rho}$ outside of $\Omega$. Passing to the limit in (\ref{912}), we see that $\rho$ is the solution to
\begin{equation}
-\int_0^T \int _{\Omega} \rho \partial_t \psi dxdt - \int_{\Omega} \rho_0 \psi(0,x )\ dx = \int_0^T \int_\Omega (\rho u) \cdot \nabla\psi \ dxdt \quad \forall \psi \in \mathcal{D}([0,T)\times \Omega). \label{619}
\end{equation}

\subsection{Induction equation}

In the induction equation, all the approximation terms already vanished during the last limit passage. Thus the limit passage with respect to $\eta \rightarrow 0$ works by the same arguments as before. Indeed, we can first use the uniform convergence (\ref{865}) to check that for any $\gamma > 0$ there exists $\delta (\gamma)>0$ such that
\begin{equation}
S_\gamma (\chi (t)) \subset S(\chi_\eta (t)) \subset S^{\frac{\gamma}{2}} (\chi (t)) \subset S^\gamma (\chi (t)) \quad \forall t \in [0,T],\ \eta < \delta (\gamma). \label{595}
\end{equation}

Then we can argue as in Section \ref{regularizationlimitinduction} to conclude
\begin{align}
\text{curl}B = 0 \quad \text{a.e. in } Q^S(\chi,T)\bigcap Q, \label{828}
\end{align}

and
\begin{align}
\left(u_\eta \times B_\eta \right) \cdot \operatorname{curl} b \rightharpoonup \left(u \times B \right) \cdot \operatorname{curl}b \quad &\text{in } L^\frac{6}{5}\left((0,T)\times \Omega\right), \label{948} \\ \operatorname{curl}B_{\eta} \times B_{\eta} \rightharpoonup \operatorname{curl} B \times B \quad &\text{in } L^\frac{6}{5}\left((0,T)\times \Omega\right) \label{929}
\end{align}

for all $b \in Y(\chi,T)$. Exploiting further the convergence (\ref{869}) of the initial data, we can pass to the limit in (\ref{914}) and obtain
\begin{align}
& - \int _0^{T} \int_\Omega B \cdot \partial_t b\ dxdt - \int_\Omega B_0 \cdot b(0,x)\ dx = \int _0^{T} \int _\Omega \left[ -\frac{1}{\sigma \mu} \text{curl} B + u \times B + \frac{1}{\sigma} J \right] \cdot \text{curl} b\ dxdt \label{596}
\end{align}

for all $b \in Y(\chi, T)$.

\subsection{Momentum equation} \label{etamomentum}

Let now $T'$ be given by (\ref{779}), i.e.\ $T'$ denotes the first time at which the rigid body $S(\chi(\cdot))$ collides with $\partial \Omega$ or, if this never happens in $[0,T]$, then $T'=T$. Since the initial distance between the body and $\partial \Omega$ is positive by (\ref{777}), the uniform convergence (\ref{865}) implies $T' > 0$ and, for any $T_0<T'$, there is some $\gamma > 0$ such that
\begin{align}
\text{dist} \Big( \partial \Omega,\ S(\chi (t))\Big) > \gamma \quad \forall t \in [0,T_0]. \label{782}
\end{align}

Our first goal in this section is to show that the limit velocity indeed coincides with a rigid velocity field in the solid region. To this end we consider an arbitrary compact set $\overline{I\times U}\subset Q^S(\chi,T')$ with an interval $I \subset (0,T')$ and some ball $U \subset \Omega$. From the first inclusion in (\ref{595}) we see that for sufficiently small $\eta$ it holds
\begin{align}
\overline{I \times U} \subset Q^S(\chi _\eta,T') \bigcap Q \quad \Leftrightarrow \quad \chi_\eta = 1 \quad \text{on } \overline{I \times U}. \nonumber
\end{align}

By the estimate (\ref{932}) this means
\begin{align}
u_\eta - \Pi_{[\chi _\eta, \rho _\eta , u_\eta]} \rightarrow 0 \quad \text{in } L^2(I \times U), \nonumber
\end{align}

and as $\overline{I \times U}$ was chosen arbitrarily we get, as desired,
\begin{equation}
u = \Pi_{[\chi, \rho, u]} \quad \text{a.e. on } Q^S(\chi,T'). \label{634}
\end{equation}

Next, we show that the projection term vanishes in the limit of the momentum equation (\ref{913}). We fix some arbitrary test function $\phi \in \mathcal{T}(\chi,T')$, i.e.\ $\phi \in \mathcal{D}([0,T)\times \Omega)$, $\operatorname{div}\phi = 0$ and there exists $\sigma > 0$ such that
\begin{align}
\mathbb{D}(\phi) = 0 \quad \text{in } \bigg\lbrace (t,x) \in (0,T')\times \Omega:\ \text{dist}\left( (t,x), \overline{Q}^S(\chi, T') \right) < \sigma \bigg\rbrace, \label{934}
\end{align}

c.f. (\ref{933}). We choose $T_0 < T'$ such that
\begin{align}
\text{supp}\phi \subset [0,T_0]\times \Omega \label{872}
\end{align}

and a corresponding $\gamma > 0$ according to (\ref{782}). By (\ref{934}) $\phi \in \mathcal{T}(\chi,T')$, there is some $0 < \sigma < \gamma$ such for all $t\in [0,T_0]$ the function $\phi(t,\cdot)$ coincides with a rigid velocity field $\phi ^S(t,\cdot)$ on $S^{\sigma}(\chi (t)) \subset \Omega$. As
\begin{align}
\chi_\eta (t,x)= 0 \quad \text{for } x \in \Omega \setminus S\left(\chi_{\eta}(t)\right), \nonumber 
\end{align}

the inclusion (\ref{595}) implies that for sufficiently small $\eta > 0$ it holds
\begin{align}
& \int_0^{T'}\int_\Omega - \frac{1}{\eta} \rho _\eta \chi _\eta \left( u_\eta - \Pi _{[\chi_\eta,\rho_\eta,u_\eta]} \right) \cdot \phi \ dxdt = \int_0^{T'}\int_\Omega - \frac{1}{\eta} \rho _\eta \chi _\eta \left( u_\eta - \Pi _{[\chi_\eta,\rho_\eta,u_\eta]} \right) \cdot \phi^S \ dxdt = 0, \label{610}
\end{align}

where the second equality is a consequence of the fact that $\Pi _{[\chi_\eta,\rho_\eta,u_\eta]}(t,\cdot)$ is the orthogonal projection of $u_\eta(t,\cdot)$ onto rigid velocity fields on $S(\chi_\eta(t))$, c.f.\ \cite[Lemma 3.1]{cottetmaitre}. \\
We further note that by the uniform bounds for $u_\eta$ in (\ref{930}) there exists a function $z_6 \in L^2(0,T;L^\frac{3}{2}(\Omega))$ such that for a chosen subsequence it holds
\begin{equation}
\rho _\eta u_\eta \otimes u_\eta \rightharpoonup z_6 \quad \text{in } L^2\left(0,T;L^\frac{3}{2}(\Omega)\right). \nonumber
\end{equation}

Combining this with the convergence (\ref{869}) of the initial data, the strong convergence (\ref{618}) of the density, the weak convergence (\ref{929}) of the magnetic term and (\ref{610}), we can pass to the limit in (\ref{913}) and obtain
\begin{align}
- \int_0^{T'} \int_\Omega \rho u \cdot \partial_t \phi \ dxdt - \int_\Omega \rho _0 u_0 \cdot \phi(0,x)\ dx =& \int _0^{T'} \int _\Omega z_6: \nabla \phi - 2\nu \mathbb{D}(u):\nabla \phi \nonumber \\
&+ \rho g\cdot \phi + \frac{1}{\mu}\left( \text{curl}B \times B \right) \cdot \phi \ dxdt \label{615}
\end{align}

for any $\phi \in \mathcal{T}(\chi,T')$. It remains to identify $z_6$. To this end it is sufficient to show that
\begin{align}
\int_0^{T_0} \int_\Omega \rho _\eta |u_\eta|^2\ dxdt \rightarrow \int_0^{T_0} \int_\Omega \rho |u|^2\ dxdt \label{fe4te4ger}
\end{align}

for arbitrary $0<T_0<T'$. Indeed, as in the proof of the classical compactness result \cite[Theorem 2.4]{lions}, this leads to strong convergence of $u_\eta$ in $L^2((0,T_0)\times \Omega)$ and in particular to
\begin{align}
z_6 = \rho u \otimes u \quad \text{a.e. on } (0,T_0)\times \Omega. \label{873}
\end{align}

Since for any arbitrary but fixed test function $\phi \in \mathcal{T}(\chi,T')$ we can find $T_0<T'$ such that the inclusion (\ref{872}) holds true, (\ref{873}) suffices to identify $z_6$ in the momentum equation (\ref{615}). The proof of (\ref{fe4te4ger}) is achieved by following mostly \cite{cottetmaitre} and using further arguments from \cite{incompressiblefeireisl}. More precisely, for fixed $0<T_0 < T'$, we choose $\gamma_\text{sup}= \gamma_\text{sup} (T_0)>0$ as the supremum over all $\gamma$ which satisfy (\ref{782}). Then for any $0 \leq \gamma \leq \frac{\gamma_\text{sup}}{4}$, $t \in [0,T_0]$ and $r \in [0,1]$ we define
\begin{align}
K_{t, \gamma}^r(\Omega) :=& \left\lbrace v(t) \in V_0^r(\Omega):\ \mathbb{D}(v(t)) = 0 \ \text{in } \mathcal{D}'\left( S^\gamma(\chi (t)) \right) \right\rbrace \label{938}
\end{align}

together with the orthogonal projection
\begin{align}
P_{\gamma}^r(t):\ H^{r}(\Omega) \rightarrow K_{t, \gamma}^r(\Omega). \label{939}
\end{align}

By the triangle inequality we estimate, for arbitrary $\psi \in \mathcal{D}(0,T_0)$, $r \in (0,1)$ and $\gamma \in (0,\frac{\gamma_{\sup}}{4}]$,
\begin{align}
&\left|\int_0^{T_0} \int_\Omega \psi \rho _\eta |u_\eta|^2\ dxdt - \int_0^{T_0} \int_\Omega  \psi \rho |u|^2\ dxdt \right| \nonumber \\
\leq& \overline{\rho} \left\| \psi \right\|_{L^\infty(0,T_0)} \left\| u_\eta \right\|_{L^2(0,T_0; L^2(\Omega)} \left\| P^r_{\gamma}u_\eta - u_\eta \right\|_{L^2(0,T_0; L^2(\Omega)} + \left| \int_0^{T_0} \int_\Omega \psi \left( \rho _\eta u_\eta \cdot P^r_{\gamma}u_\eta - \rho u \cdot P^r_{\gamma}u \right)\ dxdt \right| \nonumber \\
&+ \overline{\rho} \left\| \psi \right\|_{L^\infty(0,T_0)} \left\| u \right\|_{L^2(0,T_0; L^2(\Omega)} \left\| P^r_{\gamma}u - u \right\|_{L^2(0,T_0; L^2(\Omega)}. \label{r34z46h64th}
\end{align}

Keeping $r \in (0,1)$ and $\gamma \in (0,\frac{\gamma_\text{sup}}{4}]$ fixed, we let first $\eta$ tend to $0$. During this procedure, the second term on the right-hand side of (\ref{r34z46h64th}) vanishes, c.f.\ Lemma \ref{feireisllimit} in the Appendix. Subsequently, by letting $\gamma$ tend to $0$, also the first and the last term on the right-hand side of (\ref{r34z46h64th}) vanish, c.f.\ Lemma \ref{projectionlimit1} in the Appendix. Finally, replacing $\psi$ by a suitable sequence of cut-off functions on $[0,T_0]$, we infer the convergence (\ref{fe4te4ger}) and hence the identity (\ref{873}). 

\subsection{Proof of the main result}

Summarizing the results from Sections \ref{etabounds}--\ref{etamomentum}, we can now finish the proof of Theorem \ref{mainresult}. The regularities of $\chi$ and $\rho$ in (\ref{665}) and (\ref{635}) follow from the choice of the spaces in (\ref{867}) and (\ref{618}). As
\begin{align}
\mathbb{D}\left( \Pi_{[\chi, \rho, u]} \right) = 0, \nonumber
\end{align}

the properties of $u$ in (\ref{631}) follow from (\ref{892}) and the relation (\ref{634}) between $u$ and $\Pi_{[\chi, \rho, u]}$, while the properties of $B$ in (\ref{776}) are given by (\ref{891}) and (\ref{828}). The transport equations (\ref{666}) and (\ref{636}) were shown in (\ref{868}) and (\ref{619}), where in (\ref{868}) the function $\Pi_{[\chi,\rho,u]}$ can indeed be replaced by $u$ due to the relation (\ref{634}) between these two functions and the fact that $\chi = 0$ outside of $Q^S(\chi,T')$. The momentum equation (\ref{632}) is satisfied according to (\ref{615}), where $z_6$ was identified in (\ref{873}). The induction equation (\ref{633}) was shown to hold true in (\ref{596}). The energy inequality (\ref{638}) follows by dropping the nonnegative projection term in the energy inequality (\ref{915}) on the $\eta$-level and exploiting the weak lower semicontinuity of norms. The isometry $X(s;t,\cdot)$ is given by $X=X^{\Pi_{[\chi, \rho, u]}}$. Indeed, by (\ref{865}), $X^{\Pi_{[\chi, \rho, u]}}$ is the (pointwise) limit of a sequence of isometries and hence an isometry itself. Moreover, from the continuity of $X^{\Pi_{[\chi, \rho, u]}}$ and the fact that $X^{\Pi_{[\chi, \rho, u]}}(s;s,\cdot) = \operatorname{Id}$ it follows that $X^{\Pi_{[\chi, \rho, u]}}$ is orientation preserving. Finally, by the group property \cite[(76)]{dipernalions}, which is satisfied by $X^{\Pi_{[\chi, \rho, u]}}$ as a solution to the initial value problem (\ref{866}), it holds that
\begin{align}
S(\chi (t)) &= \left\lbrace x \in \mathbb{R}^3:\ \chi(t,x)=1 \right\rbrace = \left\lbrace X^{\Pi_{[\chi, \rho, u]}} (0;t,x):\ x \in S \right\rbrace = X^{\Pi_{[\chi, \rho, u]}} (0;t,S) \nonumber \\
&= X^{\Pi_{[\chi, \rho, u]}}\left(s;t,\left\lbrace X^{\Pi_{[\chi, \rho, u]}} (0;s,S) \right\rbrace \right) = X^{\Pi_{[\chi, \rho, u]}}(s;t,S(\chi (s))) \nonumber
\end{align}

for all $s,t \in [0,T']$. This yields the identity (\ref{780}) and thus concludes the proof.

{\centering \section{Appendix} \label{appendix} \par 
}

In the limit passage with respect to $\Delta t \rightarrow 0$ the following variant of \cite[Theorem 8.9]{roubicek} is used, which guarantees that the weak limits of different interpolants of the same discrete functions coincide. \\

\begin{satz}
\label{equalityofrothelimits}
Let $f_{\Delta t}, \overline{f}_{\Delta t}, \overline{f}_{\Delta t}'$ be piecewise affine and, respectively, piecewise constant interpolants of discrete functions $f^k_{\Delta t}$, $k=0,...,\frac{T}{\Delta t}$ defined as in (\ref{21})--(\ref{23}). Assume further that
\begin{align}
f_{\Delta t} \buildrel\ast\over\rightharpoonup f \quad \text{in } L^\infty(0,T;L^2(\Omega)),\quad \overline{f}_{\Delta t} \buildrel\ast\over\rightharpoonup \overline{f} \quad \text{in } L^\infty(0,T;L^2(\Omega))\quad \overline{f}'_{\Delta t} \buildrel\ast\over\rightharpoonup \overline{f}' \quad \text{in } L^\infty(0,T;L^2(\Omega)) \nonumber
\end{align}

Then it holds
\begin{align}
f = \overline{f} = \overline{f}'. \label{918}
\end{align}

\end{satz}

The proof of the first identity in \eqref{918} can be found in \cite[Theorem 8.9]{roubicek}, the second identity follows essentially by the same arguments. For the limit passage in the transport equation we use the following result, which is a variant of \cite[Lemma 5.2, Corollary 5.2, Corollary 5.3]{tucsnak}. A proof of the assertion is given therein.

\begin{satz}
\label{strongconvergencecharacteristics}

Assume that for any $n \in \mathbb{N}$, the function
\begin{align}
\Pi_n: [0,T] \times \mathbb{R}^3 \rightarrow \mathbb{R}^3,\quad \Pi_n(t,x) := v_n(t) + w_n(t) \times x,\quad v_n, w_n \in L^\infty(0,T), \nonumber
\end{align}

satisfies
\begin{align}
\left\| v_n \right\|_{L^\infty (0,T)},\ \left\| w_n \right\|_{L^\infty (0,T)} \leq c \label{841}
\end{align}

with $c$ independent of $n$. Denote further by $X_n$ the Carathéodory solution of
\begin{align}
\frac{dX_n(s;t,x)}{dt} = \Pi_n \left(t, X_n(s;t,x) \right), \quad \quad X_n(s;s,x) = x,\quad x \in \mathbb{R}^3, \label{845}
\end{align}

and by $\chi_n(t,x) = \chi_0(X_n(t;0,x))$ the corresponding solution to
\begin{equation}
-\int_0^T \int_{\mathbb{R}^3 } \chi_n \partial_t \Theta dxdt - \int _{\mathbb{R}^3 } \chi _0 \Theta (0,x)\ dx = \int_0^T \int_{\mathbb{R}^3} \left( \chi_n \Pi_n \right) \cdot \nabla \Theta \ dxdt \quad \forall \Theta \in \mathcal{D}([0,T)\times \mathbb{R}^3). \label{850}
\end{equation}

Then, passing to subsequences if necessary, it holds that
\begin{align}
X_n \rightarrow X \ \ \ &\text{in } C\left([0,T]\times [0,T];C_\text{loc}\left(\mathbb{R}^3 \right)\right), \label{848} \\
\chi_n \rightarrow \chi \ \ \ &\text{in } C([0,T];L_\text{loc}^p(\mathbb{R}^3 ))\quad \forall 1 \leq p < \infty \label{843}
\end{align}

with $X$ denoting the unique solution of
\begin{align}
\frac{dX(s;t,x)}{dt} = \Pi \left(t, X (s;t,x) \right), \quad \quad X(s;s,x) = x, \label{280}
\end{align}

$\chi$ the one of
\begin{equation}
-\int_0^T \int_{\mathbb{R}^3 } \chi \partial_t \Theta dxdt - \int _{\mathbb{R}^3 } \chi _0 \Theta (0,x)\ dx = \int_0^T \int_{\mathbb{R}^3} \left( \chi \Pi \right) \cdot \nabla \Theta \ dxdt \quad \forall \Theta \in \mathcal{D}([0,T)\times \mathbb{R}^3) \label{844}
\end{equation}

and with $\Pi$ given by
\begin{equation}
\Pi_n \buildrel\ast\over\rightharpoonup \Pi \ \ \ \text{in } L^\infty(0,T;W^{1,\infty}_\text{loc} (\mathbb{R}^3)),\quad \quad \Pi(t,x) = v(t) + w(t) \times x,\quad \quad v,w \in L^\infty(0,T). \label{842}
\end{equation}

Moreover,
\begin{equation}
\chi (t, x) = \chi_0\left(X(t;0,x)\right). \label{118}
\end{equation}

\end{satz}

\begin{remark}
\label{strongconvergencecharacteristicsremark}
If $\chi_0$ has compact support in $\mathbb{R}^3$, the relation $\chi_n(t,x) = \chi_0 (X_n(t;0,x))$ allows us to improve the local convergence (\ref{843}) to
\begin{align}
\chi_n \rightarrow \chi \ \ \ &\text{in } C([0,T];L^p(\mathbb{R}^3 ))\quad \forall 1 \leq p < \infty . \nonumber
\end{align}
\end{remark}

In the limit passage with respect to $\eta \rightarrow 0$, we exploit the following two results, Lemma \ref{feireisllimit} and Lemma \ref{projectionlimit1}, which are versions of \cite[Lemma 3.4]{incompressiblefeireisl} (c.f.\ \cite[Lemma 3.8]{cottetmaitre} for a related result) and \cite[Lemma 3.6, Lemma 3.7]{cottetmaitre} respectively. Detailed proofs of the results can be found in these references.

\begin{satz} \label{feireisllimit}

Let $0<T_0<T'$ be fixed, where $T'$ is defined by (\ref{779}). Let further $\gamma_\text{sup}= \gamma_\text{sup} (T_0)>0$ be the supremum of all $\gamma$ which satisfy (\ref{782}). Then, for any $\gamma \in (0,\frac{\gamma_\text{sup}}{4}]$ and any $0<r<1$, it holds
\begin{align}
\left| \int_0^{T_0} \int_\Omega \psi \left( \rho _\eta u_\eta \cdot P^r_{\gamma}u_\eta - \rho u \cdot P^r_{\gamma}u \right)\ dxdt \right| \rightarrow 0 \quad \text{for } \eta \rightarrow 0. \nonumber
\end{align}

\end{satz}

\begin{satz} \label{projectionlimit1}

Let $T_0$ be as in Lemma \ref{feireisllimit}. For any fixed $r \in (0,1)$ it holds
\begin{align}
(i)&\quad \lim_{\gamma \rightarrow 0} \lim_{\eta \rightarrow 0} \left\| P^r_{\gamma}u_\eta - u_\eta \right\|_{L^2(0,T_0;L^2(\Omega))} = 0, \nonumber \\
(ii)&\quad \lim_{\gamma \rightarrow 0} \left\| P^r_{\gamma}u - u \right\|_{L^2(0,T_0;L^2(\Omega))} = 0. \nonumber
\end{align}

\end{satz}

\begin{center}
\Large\textbf{Acknowledgements} \\[4mm]
\end{center}

\textit{This work has been supported by the Czech Science Foundation (GA\v CR) through projects 22-08633J (for \v S.N. and J.S.) and GJ19–11707Y (for B.B.). Also, the work of B.B. and J.S. was supported by the Charles University, Prague via project PRIMUS/19/SCI/01. Moreover, \it \v S. N. has been supported by  Praemium Academiæ of \v S. Ne\v casov\' a. Finally, the Institute of Mathematics, CAS is supported by RVO:67985840.}

\end{document}